\documentclass{amsart}

\usepackage{amssymb}
\usepackage{amsmath}
\usepackage{graphicx}
\usepackage{color,epsfig}
\usepackage{mathdots}
\usepackage[dvipdfm,colorlinks,backref,pagebackref,hypertexnames=false,linkcolor=blue,urlcolor=green,citecolor=red]{hyperref}

\newcommand{\Z}{{\mathbb{Z}}}
\newcommand{\F}{{\mathbb{F}}}

\newcommand{\la}{\langle}
\newcommand{\ra}{\rangle}

\newcommand{\lam}{\lambda}
\newcommand{\Sig}{\Sigma}
\newcommand{\sig}{\sigma}

\renewcommand{\epsilon}{\varepsilon}
\newcommand{\Irr}{{\rm Irr}}
\newcommand{\GL}{{\rm GL}}

\newtheorem{theorem}{Theorem}[section]
\newtheorem{lemma}[theorem]{Lemma}

\newtheorem{corollary}[theorem]{Corollary}

\newtheorem{stheorem}{Theorem}[subsection]
\newtheorem{slemma}[stheorem]{Lemma}

{\theoremstyle{remark} \newtheorem{definition}[theorem]{Definition}}

\begin{document}


\title[Counting irreducible representations of large degree]{Counting irreducible representations of large degree of the upper triangular groups}


\author{Tung Le}

\email{lttung96@yahoo.com}

\address{Institute of Mathematics, University of Aberdeen, Aberdeen AB24 3UE, U.K.}

\date{March 10, 2010}

\begin{abstract}

Let $U_n(q)$ be the upper triangular group of degree $n$ over the finite field $\F_q$ with $q$ elements. In this paper, we present constructions of large degree ordinary irreducible representations of $U_n(q)$ where $n\geq 7$, and then determine the number of irreducible representations of largest, second largest and third largest degrees.

\end{abstract}





\maketitle



\section{Introduction}
Let $q$ be a power of a prime $p$ and $\F_q$ a field with $q$
elements. The group $U_n(q)$ of all upper triangular $(n\times
n)$-matrices over $\F_q$ with all diagonal entries equal to $1$ is a
Sylow $p$-subgroup of $\GL_n(\F_q)$. It was conjectured by Higman
\cite{Hig} that the number of conjugacy classes of $U_n(q)$ is given
by a polynomial in $q$ with integer coefficients. Isaacs \cite{Is2} showed that the degrees of all
irreducible characters of $U_n(q)$ are powers of $q$. Huppert \cite{Hup} proved that character degrees of $U_n(q)$ are precisely of the form $\{q^e:0\leq e \leq \mu(n)\}$ where the upper bound $\mu(n)$ was known to Lehrer \cite{Leh}.
Lehrer \cite{Leh} conjectured that each number  $N_{n,e}(q)$ of
irreducible characters of $U_n(q)$ of degree $q^e$ is given by a
polynomial in $q$ with integer coefficients. Isaacs \cite{Is3}
suggested a strengthened form of Lehrer's conjecture stating that
$N_{n,e}(q)$ is given by a polynomial in $(q-1)$ with nonnegative
integer coefficients. So, Isaacs' conjecture implies Lehrer's Conjecture, which in turn implies Higman's Conjecture.

To study irreducible characters of $U_n(q)$, Andr\'e
\cite{Carlos2} constructs elementary characters using the Kirillov
orbit method, which is similar to Lehrer's construction \cite{Leh}. He then defines basic characters as tensor products of sets of
elementary characters. The set nonprincipal irreducible characters $\Irr(U_n(q))^*$ of
$U_n(q)$ is partitioned by the irreducible constituents of basic characters. Diaconis and Isaacs \cite{Is-Di} generalize the notion of basic characters to that of supercharacters.

When $q$ is a proper prime power, $U_n(q)$ has no faithful
characters. We define { almost faithful}
characters as follows.

\begin{definition}
\label{def:afaithful}
Let $G$ be a group with center $Z(G)$. A character $\chi$ of $G$
is called {\em almost faithful} if $\chi|_{Z(G)}\neq \chi(1) 1_{Z(G)}$.
\end{definition}

It is clear that when $q=p$, $Z(U_n(q))$ is cyclic, almost faithful implies faithful. Here, we
approach complex representations of a Sylow $p$-subgroup $U$ of a finite
group of Lie type $G(q)$ by its root system. For each positive root,
we construct a minimal degree almost faithful irreducible character,
which is called {\it midafi}. The constituents of a certain collection of tensor
products of midafis partition $\Irr(U)^*$. If $G=\GL_n(q)$, it turns out that our midafi characters are precisely the elementary characters of Andr\'e, and moreover, the resulting partitions of $\Irr(U_n(q))^*$ coincide. This is no longer true when $G$ is not of type $A$, see \cite{Carlos-Neto,Frank-Tung_Kay}.

To compute the numbers $N_{n,e}(q)$ of $U_n(q)$, Isaacs
\cite{Is3} uses combinatorial methods to
count the numbers of orbits in a certain nilpotent associative $\F_q$-algebra. He obtains the exact numbers of largest and second largest degree irreducible characters of $U_n(q)$. Here, we decompose certain
supercharacters to construct all irreducible characters of $U_n(q)$ of degree $q^a$, $a\in \{\mu(n),\mu(n)-1,\mu(n)-2\}$. We prove the following.

\begin{theorem} \label{intro:main-thm}
The following formulae hold.
\begin{itemize}
\item[(i)] For $n\geq 3$, we have \[N_{n,\mu(n)}(q)=(q-1) N_{n-2,\mu(n-2)}(q).\]

\item[(ii)] For $n\geq 5$, we have \[N_{n,\mu(n)-1}(q)=(q-1)N_{n-2,\mu(n-2)-1}(q)+q(q-1)^2N_{n-4,\mu(n-4)}(q).\]

\item[(iii)] For $n\geq 7$, we have
\[\begin{array}{ll}
N_{n,\mu(n)-2}(q)=&(q-1)N_{n-2,\mu(n-2)-2}(q)+q(q-1)^2N_{n-4,\mu(n-4)-1}(q)\\
&+(2q^2+q-1)(q-1)^3N_{n-6,\mu(n-6)}(q).
\end{array}
\]
\end{itemize}
\end{theorem}

Using the initial values from \cite{Is3}, we show that the number $N_{n,\mu(n)-2}(q)$ of third largest degree irreducible characters of $U_n(q)$ is a polynomial in $(q-1)$ with positive integer coefficients.

In \cite{Mar}, Marjoram uses combinatorial methods to establish the exact formula for $N_{n,\mu(n)}(q)$. Hence, Theorem \ref{intro:main-thm} (i) is properly included in Marjoram's and Isaacs' work. But, since we prove part (i) directly from the recursion construction of $U_n(q)$ to $U_{n-2}(q)$ by using supercharacters, it seems appropriate to present it here. In addition, Marjoram's thesis includes polynomial formulas for $N_{n,\mu(n)-1}(q)$ and $N_{n,\mu(n)-2}(q)$ when $n$ is even.

At the end of Section \ref{section:three_largest}, we give a very sketch of how our methods can be pushed further to compute the number $N_{n,\mu(n)-k}(q)$ of $U_n(q)$ for $n\geq 2k+1$. However, the length of the full argument exceeds the scope of the present paper. Without proof, we state that
\[\begin{array}{ll}
N_{n,\mu(n)-3}(q)=&(q-1)N_{n-2,\mu(n-2)-3}(q)+q(q-1)^2N_{n-4,\mu(n-4)-2}(q)\\
& +(2q^2+q-1)(q-1)^3N_{n-6,\mu(n-6)-1}(q)+q^2(q-1)^3 N_{n-6,\mu(n-6)}(q)\\
&+(4q^3+4q^2-3q-1)(q-1)^4N_{n-8,\mu(n-8)}(q).
\end{array}
\]
The proof of Theorem \ref{thm:Andre} shows how to decompose a tensor product
of two supercharacters into a sum of supercharacters. Notice that Theorem \ref{thm:Andre} is similar to \cite[Theorem 1]{Carlos2}.

This paper is organized as follows.  We first introduce the relevant properties of
$U_n(q)$ and supercharacters. Next, we prove Theorem \ref{thm:Andre} by using our language and constructions. In Section \ref{section:three_largest}, we construct and count the large degree characters.


\section{Some properties of $U_n(q)$'s characters} \label{section:summary}

Let $G$ be a group. Denote $G^*:=G-\{1\}$, $\Irr(G)$ the set of all complex irreducible characters of
$G$, and $\Irr(G)^*:=\Irr(G)-\{1_G\}$. For $H\unlhd G$, let $\Irr(G/H)$ denote the irreducible character set of $G$ with $H$ in the kernel. If $K\leq G$ such that $G=H\rtimes K$, then for each character $\xi$ of $K$, we denote the inflation of $\xi$ to $G$ by $\xi_G$, i.e. $\xi_G$ is the extension of $\xi$ to $G$  with $H\subset \ker(\xi_G)$. For $H\leq G$ and $\xi\in \Irr(H)$, define $\Irr(G,\xi):=\{\chi\in \Irr(G):(\chi,\xi^G)\neq 0\}$. Furthermore, for a character $\chi$ of $G$, we denote its restriction to $H$ by $\chi|_H$.

Let $\Sig:=\Sig_{n-1}=\la \alpha_1,...,\alpha_{n-1}\ra$ be the root system
of $\GL_n(q)$ with respect to the maximal split torus equal to the
diagonal group, see \cite{cart2}. Denote $\alpha_{i,j}:=\alpha_i+\alpha_{i+1}+...+\alpha_j$ for all $0<i\leq j<n$, and define the {\it height} of $\alpha_{i,j}$ as $ht(\alpha):=j-i+1$. Let $\Sig^+$ denote the set of all positive roots and
$X_\alpha$ the root subgroup of $\GL_n(q)$ corresponding to $\alpha\in \Sig^+$. Hence,
$X_{\alpha_{i,j}}$ is the set of all upper triangular matrices in the
form $I_n+c\cdot e_{i,j+1}$, where $I_n$ is the  identity $(n\times
n)$-matrix, $c\in \F_q$ and $e_{i,j+1}$ is the zero matrix except a `$1$' at entry
$(i,j+1)$. The {\it upper triangular} group $U_n(q)$ is generated by
all $X_\alpha$'s, $\alpha\in\Sig^+$. We write $U$ for $U_n(q)$ if $n$ and $q$ are clear from the context.

For each $\alpha\in\Sig^+$, there is a corresponding set of $(q-1)$ orbits of the nilpotent Lie algebra of $U$. For each such orbit,
 Andr\'e \cite{Carlos2} constructs an irreducible character, which he calls {\it elementary.} We shall reconstruct these characters by using the root system.

\begin{definition}
\label{def:fundamental}
For each $\alpha\in\Sig^+$, we define
\begin{itemize}

\item[(i)] $arm(\alpha):=\{\beta\in\Sig^+:\alpha,\beta\mbox{ on the same row and } (\alpha-\beta)\in\Sig^+\}$ and\\ $A_\alpha:=\la X_\beta:\beta\in arm(\alpha)\ra$,

\item[(ii)] $leg(\alpha):=\{\beta\in\Sig^+:\alpha,\beta\mbox{ on the same column and } (\alpha-\beta)\in\Sig^+\}$ and\\ $L_\alpha:=\la X_\beta:\beta\in leg(\alpha)\ra$,

\item[(iii)]  the hook $h(\alpha):=\{\alpha\}\cup arm(\alpha) \cup leg(\alpha)$, and  the hook group $H_\alpha:=\la X_\beta:\beta\in h(\alpha)\ra$,

\item[(iv)] the base group $V_\alpha:=\la X_\beta:\beta\in\Sig^+- arm(\alpha)\ra$,

\item[(v)] the subtriangular group $U_\alpha:=\la X_\beta: \beta\in\Sig^+,\exists \gamma\in\Sig^+\cup\{0\}, (\beta+\gamma)\in h(\alpha)\ra$,

\item[(vi)]  the radical $R_\alpha:=\la X_\beta: \beta\in\Sig^+, X_\beta\nsubseteq U_\alpha \ra$.
\end{itemize}

\end{definition}

If $\alpha$ is a simple root, $arm(\alpha)$, $leg(\alpha)$ are empty and $A_\alpha$, $L_\alpha$ are trivial. If $ht(\alpha)>1$, then $A_\alpha$, $L_\alpha$ are abelian, $[A_\alpha,L_\alpha]=X_\alpha=Z(H_\alpha)$ and $H_\alpha=A_\alpha X_\alpha\rtimes L_\alpha=A_\alpha\ltimes X_\alpha L_\alpha$ is a special group. We have $U=U_{\alpha_{1,n-1}}=R_\alpha\rtimes U_\alpha$ for all $\alpha\in\Sig^+$. Therefore, all characters of $U_\alpha$ are inflated to $U$. These groups can be visualized as follows.
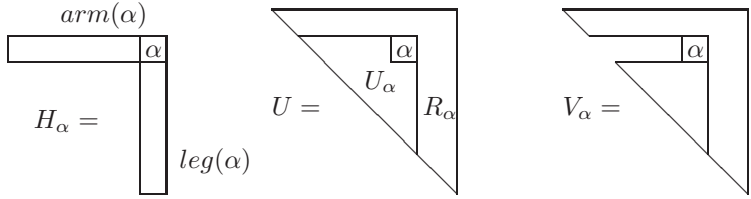
\begin{center}  
\begin{picture}(270,80)

\put(0,60){\line(1,0){60}}
\put(60,60){\line(0,-1){60}}
\put(60,0){\line(-1,0){10}}
\put(50,0){\line(0,1){60}}
\put(0,60){\line(0,-1){10}}
\put(0,50){\line(1,0){60}}
\put(52,52){$\alpha$}
\put(10,25){$H_\alpha=$}
\put(20,65){$arm(\alpha)$}
\put(65,10){$leg(\alpha)$}

\put(100,30){$U=$}
\put(100,70){\line(1,0){70}}
\put(170,70){\line(0,-1){70}}
\put(170,0){\line(-1,1){70}}
\put(155,60){\line(0,-1){45}}
\put(155,60){\line(-1,0){45}}
\put(145,60){\line(0,-1){10}}
\put(155,50){\line(-1,0){10}}

\put(147,52){$\alpha$}
\put(135,40){$U_\alpha$}
\put(157,30){$R_\alpha$}

\put(210,30){$V_\alpha=$}
\put(210,70){\line(1,0){70}}
\put(280,70){\line(0,-1){70}}
\put(280,0){\line(-1,1){50}}
\put(265,60){\line(0,-1){45}}
\put(265,60){\line(-1,0){45}}
\put(255,60){\line(0,-1){10}}
\put(265,50){\line(-1,0){35}}
\put(210,70){\line(1,-1){10}}
\put(257,52){$\alpha$}

\end{picture}
\end{center}

Since $X_\alpha\cap [V_\alpha,V_\alpha]=\{1\}$, there exist
$(q-1)$ linear characters $\lam\in \Irr(V_\alpha/[V_\alpha,V_\alpha])$ such
that $\lam|_{X_\alpha}\neq 1_{X_\alpha}$ and $\lam|_{X_\beta}=1_{X_\beta}$
for all $X_{\beta}\subset V_\alpha$ with $\beta\neq\alpha$. From now on, we denote by
$\Irr(V_\alpha/[V_\alpha,V_\alpha])^*$ the set of these $(q-1)$ linear characters of
$V_\alpha$. For each $\lam\in \Irr(V_\alpha/[V_\alpha,V_\alpha])^*$, we call $\lam^U$ an {\it elementary character} associated to $\alpha$. Since the hook group $H_\alpha$ is special, the following lemma is obvious. However, we provide a short proof.

\begin{lemma} \label{lem:midafi}

$\lam^U$ is irreducible and $((\lam|_{V_\alpha\cap
U_\alpha})^{U_\alpha})_U=\lam^U$. Moreover, $\lam^U|_{H_\alpha}$
is irreducible and extendable to $U$.

\end{lemma}

\begin{proof}
If $\alpha$ is a simple root since $H_\alpha=U_\alpha=X_\alpha$, $V_\alpha=U$ and $\lam^U=\lam$. Now we suppose that the height $ht(\alpha)>1$. Since  $R_\alpha\subset \ker(\lam^U)$, it is enough to show that $\lam^U|_{H_\alpha}$
is irreducible. By Mackey formula for the double coset $U=V_\alpha
H_\alpha$, $\lam^U|_{H_\alpha}={\lam|_{V_\alpha\cap H_\alpha}}^{H_\alpha}$.

Here $V_{\alpha}\cap H_\alpha=X_{\alpha}L_\alpha$, and $H_\alpha=X_\alpha L_\alpha\rtimes A_\alpha$. Each $x\in A_\alpha^*$ can be written uniquely as
$x_{\beta_1}...x_{\beta_k}$ where $\beta_i\in arm(\alpha)$,
$ht(\beta_i)<ht(\beta_{i+1})$ and $x_{\beta_i}\in X_{\beta_i}^*$. Let $\gamma\in
leg(\alpha)$ such that $\gamma+\beta_1=\alpha$. For all $y\in
X_\gamma$, we have $[y,x]\in X_{\alpha}$ and
\[({}^x\lam|_{X_\alpha L_\alpha})(y)=\lam(y^x)=\lam(y[y,x])=\lam([y,x])\lam(y)\neq\lam|_{X_\alpha L_\alpha}(y)=1.\]
So the inertia group
$I_{H_\alpha}(\lam|_{X_\alpha L_\alpha})=X_\alpha L_\alpha$. By Clifford theory,
${\lam|_{X_\alpha L_\alpha}}^{H_\alpha}$ is
irreducible.
\end{proof}

From the proof of Lemma \ref{lem:midafi}, the following corollary is clear.

\begin{corollary} \label{cor:midafi_to_leg}
For each $\alpha\in \Sig^+$, the following hold.
\begin{itemize}

\item[(i)] $H_\alpha$ has $q^{2\cdot|arm(\alpha)|}$ linear characters and $(q-1)$ almost faithful irreducible characters of degree $q^{|arm(\alpha)|}=|A_\alpha|=|L_\alpha|$, which are $\{\lam^U|_{H_\alpha}:\ \lam\in \Irr(V_\alpha/[V_\alpha,V_\alpha])^*\}$.

\item[(ii)] $\lam^U|_{L_\alpha}={1_{\{1\}}}^{L_\alpha}$ the regular character of $L_\alpha$.

\item[(iii)] $\varphi^{H_\alpha} =\lam^U|_{H_\alpha}$ for all
$\varphi\in \Irr(X_{\alpha}L_\alpha,\lam|_{X_\alpha})$.
\end{itemize}
\end{corollary}

Two roots $\alpha$, $\beta\in\Sig^+$ are called {\it separate} if they
are on neither the same row nor same column. A nonempty subset $D$
of $\Sig^+$ is called a {\it basic set} if all roots in $D$ are pairwise separate. For each basic set $D$, let $E_D=\bigoplus_{\alpha\in D} \Irr(V_\alpha/[V_\alpha,V_\alpha])^*$. For each basic set $D$ and $\phi\in
E_D$, Andr\'e \cite{Carlos2} defines a {\it basic character} as
\[\xi_{D,\phi}=\bigotimes_{\lam_\alpha\in\phi}{\lam_\alpha}^U.\]

In Diaconis and Isaacs \cite{Is-Di}, $\xi_{D,\phi}$ is an example of a {\it supercharacter}.
It is known that all $p$-groups are monomial, i.e. all
irreducible characters are induced from linear characters. Here, supercharacters may not be irreducible, but the following lemma shows that they are
always induced from linear characters.

\begin{definition}
\label{def:Vb+Ld}
For each basic set $D$ and $\phi\in E_D$, we define \[V_D:=\bigcap_{\alpha\in D}
V_\alpha\ \  \mbox{ and }\ \ \lam_D:=\bigotimes_{\lam_\alpha\in \phi}\lam_\alpha|_{V_D}.\]

\end{definition}

\begin{lemma} \label{lem:tobelinear}

We have
$\xi_{D,\phi}=\bigotimes_{\lam_\alpha\in \phi}{\lam_\alpha}^U={\lam_D}^U$.

\end{lemma}

\begin{proof}
We induct on the size of $D$. Suppose that
$D':=D\cup\{\beta\}$ is a basic set and for all $\phi\in E_D$,
$\xi_{D,\phi}={\lam_D}^U$ where $\lam_D=\bigotimes_{\lam\in \phi}
\lam|_{V_D}$ and $V_D=\bigcap_{\alpha\in D}V_\alpha$. Let $\phi':=\phi\cup
\{\lam_\beta\}$ for some $\phi\in E_D$ and $\lam_\beta\in
\Irr(V_\beta/[V_\beta,V_\beta])^*$. Since $D\cup\{\beta\}$ is a basic set, $U=V_D V_\beta$, by Mackey formula for this
double coset, we have
\[\begin{array}{lll}
 \xi_{D',\phi'}&=&\xi_{D,\phi}\otimes \lam_{\beta}^U\\
  &=&\lam_D^U\otimes \lam_\beta^U\\
 &=&(\lam_D^U|_{V_\beta} \otimes \lam_\beta)^U\\
 &=& ({\lam_D|_{V_D\cap V_\beta}}^{V_\beta}\otimes\lam_\beta )^U\\
 &=&((\lam_D|_{V_D\cap V_\beta}\otimes \lam_\beta|_{V_D\cap V_\beta})^{V_\beta})^U\\
 &=&{\lam_{D'}}^U,
\end{array}\]
where $\lam_{D'}=\bigotimes_{\lam\in \phi'}\lam|_{V_{D'}}$ and $V_{D'}=\bigcap_{\alpha\in D'}V_\alpha$.
\end{proof}

Actually in the proof of Lemma \ref{lem:tobelinear}, we only need `no two roots in $D$ on the same row'. The main motivation for defining supercharacters is to get a partition of the set of all nonprincipal irreducible characters $\Irr(U)^*$ of $U$.

\begin{theorem}
\label{thm:Andre}
For each $\chi\in \Irr(U)^*$, there exist uniquely  a basic set $D$ and a set $\phi\in E_D$ such that $\chi$ is a constituent of $\xi_{D,\phi}$.

\end{theorem}

Notice that this is \cite[Theorem 1]{Carlos2}. In Section \ref{section:supercharacter}, we reprove it by using different methods. In Section \ref{section:three_largest}, our methods are extended to prove Theorem \ref{intro:main-thm}.


\section{Supercharacters} \label{section:supercharacter}

When $q$ is not prime, $Z(U)$ is not cyclic, $U$ has no faithful characters. That explains why we come up with `{\it almost faithful}' in Definition \ref{def:afaithful}. Notice that for a group $G$, $Z(G)\subset Z(\chi)$ for all $\chi\in \Irr(G)$, where $Z(\chi):=\{g\in G:|\chi(g)|=\chi(1)\}$. It is clear that ${\lam_\alpha}^U|_{U_\alpha}$ is a minimal degree almost
faithful irreducible character of $U_\alpha$ for all $\lam_\alpha\in
\Irr(V_\alpha/[V_\alpha,V_\alpha])^*$. Recall that if $\chi\in \Irr(G)$ is almost faithful, there exists uniquely $\sig\in
\Irr(Z(G))^*$ such that $\chi\in \Irr(G,\sig)$. The next theorem characterizes all
almost faithful irreducible characters.

\begin{theorem}
\label{thm:unique_factor} 
Each almost faithful irreducible character of $U$   factors uniquely
into a tensor product of an almost faithful elementary character and
a character of $\Irr(U/H_{\alpha_{1,n-1}})$.
\end{theorem}

\begin{proof}
It is known that $Z(U)=X_{\alpha_{1,n-1}}$. Let $\sig\in \Irr(X_{\alpha_{1,n-1}})^*$ and
 $\lam^U\in \Irr(U,\sig)$ be an elementary character associated to $\alpha_{1,n-1}$. Denote $H:=H_{\alpha_{1,n-1}}$. We have
\[(\lam^U|_{H})^U=(\lam^U|_{H}\otimes 1_{H})^U=\lam^U\otimes {1_{H}}^U=\lam^U\otimes\sum_{\eta\in \Irr(U/H)}\eta(1)\, \eta.\]
%

By Lemma \ref{lem:midafi}, $\lam^U|_{H}$ is irreducible and extendable to
$U$, by Clifford theory,  $\lam^U\otimes \eta$ are
irreducible for all $\eta\in \Irr(U/H)$. Since
$\sig^{H}=\lam^U(1)\, \lam^U|_{H}$, by the transitivity of
induction, $\sig^U=(\sig^{H})^U=\sum_{\eta\in
\Irr(U/H)}\lam^U(1)\eta(1)\,\lam^U\otimes \eta$. Therefore, for each
almost faithful irreducible character $\chi\in \Irr(U,\sig)$, there
exists unique $\eta\in \Irr(U/H)$ such that $\chi=\lam^U\otimes \eta$.
\end{proof}

For each appropriate number $e\geq n-2$, by Theorem \ref{thm:unique_factor}, the set of all irreducible characters of degree $q^e$ always contains almost faithful ones. We will apply Theorem \ref{thm:unique_factor} to construct the representations of large degree irreducible characters in Section \ref{section:three_largest}.

If $\chi\in\Irr(U)$ is not almost faithful,  $X_{\alpha_{1,n-1}}\leq\ker(\chi)$. So we can work with $U/X_{\alpha_{1,n-1}}$ when decomposing supercharacters $\xi_{D,\phi}$ where $\alpha_{1,n-1}\notin D$. First, we study tensor products of two elementary characters. A tensor product of an elementary character with a linear character clearly remains irreducible. In particular, if ${\lam_\alpha}^U$, ${\lam_\beta}^U$ are elementary characters where $\alpha\in arm(\beta)$ (or $leg(\beta)$) and $\alpha$ is a simple root, i.e. ${\lam_\alpha}^U=\lam_\alpha$, then 
$\lam_\alpha\otimes{\lam_\beta}^U={\lam_\beta}^U$.

\begin{lemma} \label{lem:prod2midafis}

Let ${\lam_1}^U$, ${\lam_2}^U$ be two nonlinear elementary
characters associated to $\alpha_{i,j}$ and $\alpha_{l,k}$ respectively. The decomposition of ${\lam_1}^U\otimes{\lam_2}^U$ is as follows:
\begin{itemize}

\item[(i)] if $|H_{\alpha_{i,j}}\cap H_{\alpha_{l,k}}|=1$ then ${\lam_1}^U\otimes {\lam_2}^U$ is irreducible,

\item[(ii)]   if $|H_{\alpha_{i,j}}\cap H_{\alpha_{l,k}}|=q$ then ${\lam_1}^U\otimes {\lam_2}^U$ has $q$ distinct irreducible constituents of degree $q^{(j-i)+(k-l)-1}$ with multiplicity $1$,

\item[(iii)]   if $\alpha_{i,j}\in leg(\alpha_{l,k})$ (or $arm(\alpha_{l,k})$), then ${\lam_1}^U\otimes {\lam_2}^U$ decomposes into ${\lam_2}^U$ and a sum of $\sum_{\lam_{\alpha}\in \Irr(V_\alpha/[V_\alpha,V_\alpha])^*}{\lam_2}^U\otimes{\lam_\alpha}^U$ for all $\alpha\in arm(\alpha_{i,j})$ (or $leg(\alpha_{i,j})$ respectively); each constituent has multiplicity $1$,

\item[(iv)]   if $\alpha_{i,j}=\alpha_{l,k}$, and ${\lam_1}^U\neq \overline{{\lam_2}^U}$, then ${\lam_1}^U\otimes {\lam_2}^U$ decomposes into 
    \[(q+(j-i-1)(q-1))(\lam_1\otimes \lam_2)^U+\sum^{U_\alpha\subset U_{\alpha_{i+1,j-1}}}_{\lam_\alpha\in \Irr(V_\alpha/[V_\alpha,V_\alpha])^*}(\lam_1\otimes \lam_2)^U\otimes{\lam_\alpha}^U.\]

\item[(v)]   if $\alpha_{i,j}=\alpha_{l,k}$, and ${\lam_1}^U=\overline{{\lam_2}^U}$, then ${\lam_1}^U\otimes {\lam_2}^U$ decomposes into
    \[1_U+\sum^{\beta_1\in arm(\alpha_{i,j}),\ \lam_{\beta_1}\in \Irr(V_{\beta_1}/[V_{\beta_1},V_{\beta_1}])^*}_{\beta_2\in leg(\alpha_{i,j}),\ \lam_{\beta_2}\in \Irr(V_{\beta_2}/[V_{\beta_2},V_{\beta_2}])^*}{\lam_{\beta_1}}^U+{\lam_{\beta_2}}^U+{\lam_{\beta_1}}^U\otimes{\lam_{\beta_2}}^U.\] 
\end{itemize}

\end{lemma}

\begin{proof}
(i) If $U_{\alpha_{i,j}}\subset U_{\alpha_{l,k}}$(
or $U_{\alpha_{l,k}}\subset U_{\alpha_{i,j}}$), it is clear by Theorem
\ref{thm:unique_factor}; if
$U_{\alpha_{i,j}}\cap U_{\alpha_{l,k}}=\{1\}$, it follows by $\Irr(U_{\alpha_{i,j}}\times U_{\alpha_{l,k}})\simeq \Irr(U_{\alpha_{i,j}})\times \Irr(U_{\alpha_{l,k}})$.

\begin{center} 
\begin{picture}(340,80)

\put(0,25){$V_1=$}
\put(40,70){\line(1,0){30}}    
\put(70,70){\line(0,-1){70}}  
\put(70,0){\line(-1,1){60}}   
\put(40,70){\line(0,-1){40}}
\put(50,70){\line(0,-1){50}}
\put(50,60){\line(-1,0){40}}
\put(38,74){$\alpha_{i,j}$}
\put(70,45){\line(-1,0){10}} \put(60,45){\line(0,-1){10}} \put(60,35){\line(1,0){10}}

\put(90,25){$V_2=$}
\put(90,70){\line(1,0){70}}  
\put(160,70){\line(0,-1){70}}  
\put(160,0){\line(-1,1){35}}  
\put(90,70){\line(1,-1){25}}  
\put(160,45){\line(-1,0){45}}
\put(160,35){\line(-1,0){35}}
\put(150,45){\line(0,-1){35}}
\put(163,37){$\alpha_{l,k}$}
\put(130,70){\line(0,-1){10}} \put(130,60){\line(1,0){10}} \put(140,60){\line(0,1){10}}

\put(183,25){$V_{12}=$}
\put(220,70){\line(0,-1){25}}
\put(230,70){\line(0,-1){25}}
\put(220,35){\line(0,-1){5}}
\put(230,35){\line(0,-1){15}}
\put(220,70){\line(1,0){30}} 
\put(230,60){\line(-1,0){40}}
\put(250,70){\line(0,-1){70}} 
\put(240,45){\line(0,-1){35}} 
\put(250,45){\line(-1,0){45}}
\put(250,35){\line(-1,0){35}}
\put(250,0){\line(-1,1){35}}   
\put(190,60){\line(1,-1){15}}  

\put(273,25){$V=$}
\put(310,70){\line(0,-1){25}}
\put(320,70){\line(0,-1){25}}
\put(310,35){\line(0,-1){5}}
\put(320,35){\line(0,-1){15}}
\put(310,70){\line(1,0){30}} 
\put(320,60){\line(-1,0){40}}
\put(340,70){\line(0,-1){70}} 
\put(330,45){\line(0,-1){35}} 
\put(340,45){\line(-1,0){45}}
\put(340,35){\line(-1,0){35}}
\put(340,0){\line(-1,1){35}}   
\put(280,60){\line(1,-1){15}}  

\put(305,70){\line(0,-1){10}} \put(305,70){\line(-1,0){10}} \put(295,70){\line(0,-1){10}}
\put(290,75){$\alpha_{i,l-1}$}

\end{picture}

\end{center}

(ii)  Since $|H_{\alpha_{i,j}}\cap H_{\alpha_{l,k}}|=q$, $\alpha_{i,j}$ and
$\alpha_{l,k}$ are separate. We suppose that $i<l<j<k$. It suffices to work with the quotient $U/(R_{\alpha_{i,j}}\cap R_{\alpha_{l,k}})$. Let $V_1,
V_2$ be the base groups of $\alpha_{i,j}$, $\alpha_{l,k}$ respectively. By
Lemma \ref{lem:tobelinear},
$\lam_1^U\otimes\lam_2^U=\lam_{12}^U$ where
$\lam_{12}=\lam_1|_{V_{12}}\otimes \lam_2|_{V_{12}}$, and
$V_{12}=V_1\cap V_2$. Since
$[V_{12},X_{\alpha_{i,l-1}}]=\{1\}$, we let $V:=\la
V_{12},X_{\alpha_{i,l-1}}\ra=V_{12}X_{\alpha_{i,l-1}}\simeq V_{12}\times
X_{\alpha_{i,l-1}}$. So $\lam_{12}^V={\lam_{12}}_V\otimes
{1_{V_{12}}}^V={\lam_{12}}_V\otimes \sum_{\rho\in
\Irr(V/V_{12})}\rho=\sum_{\rho\in
\Irr(X_{\alpha_{i,l-1}})}{\lam_{12}}_V\otimes\rho_V$. Therefore,
$\lam_{12}$ extends to $q$ distinct linear characters of $V$. Let
$\lam$ be an extension of $\lam_{12}$ to $V$.

Set $M_0:=V$, $M_t:=M_{t-1}X_{\alpha_{i,j-t}}$, for all $t\in[1,j-i]$, and $M_t:=M_{t-1}X_{\alpha_{l,k+(j-i)-t}}$ for all $t\in[(j-i)+1,(j-i)+(k-l)]$. Hence, $M_{l-j}=M_{l-j+1}$,
$M_{j-i}=V_2$, and $M_{(j-i)+(k-l)}=U$. To show that $\lam^U$ is
irreducible, by the transitivity of inductions, we proceed
$(j-i)+(k-l)-1$ steps of inductions along the arms of $\alpha_{i,j}$
and $\alpha_{l,k}$ respectively namely from $M_0$ to $M_1$, $M_2$, ...,
$M_{(l-k)+(j-i)}$.

Let $K_0{:=}(leg(\alpha_{i,j})\cup leg(\alpha_{l,k}))-
\{\alpha_{l,j}\}$ and $K_t{:=}K_{t-1}-\{\alpha_{j-t+1,j}\}$ for all
$t\in[1,j-i]$, $K_t:=K_{t-1}- \{\alpha_{k+(j-i)-t+1,l}\}$
for all $t\in[(j-i)+1,(j-i)+(k-l)]$. Hence,
$K_{j-l}=K_{j-l+1}$, $K_{j-i}=leg(\alpha_{l,k})$, $K_{(j-i)+(k-l)}$ is
empty, $X_\beta\subset \ker(\lam)$ for all $\beta\in K_0$. For
 $M_t=M_{t-1}X_\alpha$ such that $M_t\neq M_{t-1}$, there exists uniquely $\beta\in K_t- K_{t-1}$ satisfying $\alpha+\beta=\alpha_{i,j}$ if
$t\leq j-i$, $\alpha+\beta=\alpha_{l,k}$ if $t\geq j-i+1$.

Suppose that $\lam^L\in \Irr(L)$ for some $L=M_t$, and
$X_\beta\subset \ker(\lam^L)$ for all $\beta\in K=K_t$. If
$t=(l-k)+(j-i)$, it is done. Otherwise, the next induction step is
from $L$ to $LX_\alpha=M_{t+1}$ where $\alpha\in arm(\tau)$, $\tau=\alpha_{i,j}$
if $t\leq j-i$ and $\tau=\alpha_{l,k}$ if $t\geq (j-i)+1$. Here, $X_\alpha$ normalizes $L$ and there exists $\beta\in K$ in the leg of $\tau$ such that
$\alpha+\beta=\tau$, $X_\beta\subset \ker(\lam^L)$. So for each $x\in X_\alpha^*$, all $x_\beta\in X_\beta$, we have
\begin{center}
$ {}^x(\lam^L)(x_\beta)=\lam^L({x_\beta}^x)=\lam^L(x_\beta[x_\beta,x])=\lam([x_\beta,x])\lam^L(x_\beta)\neq \lam^L(x_\beta)=\lam^L(1)$,
\end{center}
since $[x_\beta,x]\in X_\tau\subset Z(\lam^L)$ and there exists some $x_\beta\in X_\beta$ such that $\lam([x_\beta,x])\neq 1$.

Hence, for all $x\in X_\alpha^*$, ${}^x(\lam^L)\neq \lam^L$ since
$X_\beta\nsubseteq \ker({}^x(\lam^L))$. By Mackey
formula for the double coset $L\backslash LX_\alpha/L$ represented by
$X_\alpha$, we have
\[(\lam^{LX_\alpha},\lam^{LX_\alpha})=\sum_{x\in X_\alpha}({}^x(\lam^L),\lam^L)=(\lam^L,\lam^L)=1, \mbox{ i.e. } \lam^{LX_\alpha}\in \Irr(LX_\alpha,\lam).\]
It is clear that $X_\gamma\subset \ker(\lam^{LX_{\alpha}})$ for all $\gamma\in K-\{\beta\}=K_{t+1}$.
Therefore, by induction on $t$, $\lam^U \in \Irr(U,\lam)$ of degree $q^{(l-k)+(j-i)-1}$ .
Now, we show that for any  two distinct extensions
$\lam$, $\lam'$ of $\lam_{12}$ to $V$, if $\lam^L\neq \lam'^L$  then
$\lam^{LX_\alpha}\neq \lam'^{LX_\alpha}$. Since $X_\beta\nsubseteq
\ker({}^x(\lam^L))$ for all $x\in X_\alpha^*$,
\[(\lam^{LX_\alpha},\lam'^{LX_\alpha})=\sum_{x\in X_\alpha}({}^x(\lam^L),\lam'^L)=(\lam^L,\lam'^L)=0.\]

(iii) Suppose $\alpha_{i,j}\in leg(\alpha_{l,k})$, hence, $i>l\geq 1$. Let $\tau_1:=\alpha_{i,j}$,
$\tau_2:=\alpha_{l,k}$, $V_1:=V_{\tau_1}$ and $V_2:=V_{\tau_2}$. By
Mackey formula for $V_2 V_1=U$, $\lam_1^U\otimes
\lam_2^U=(\lam_1\otimes \lam_2^U|_{V_1})^U=(\lam_1\otimes
{\lam_2|_{V_1\cap V_2}}^{V_1})^U$. By Corollary
\ref{cor:midafi_to_leg}, $\lam_1\otimes {\lam_2|_{V_1\cap
V_2}}^{V_1}=(\lam_1|_{V_1\cap V_2}\otimes \lam_2|_{V_1\cap
V_2})^{V_1}=(\lam_2|_{V_1\cap V_2})^{V_1}=\lam_2^U|_{V_1}$.
So
$\lam_1^U\otimes\lam_2^U=(\lam_2^U|_{V_1})^U=(\lam_2^U|_{V_1}\otimes
1_{V_1})^U=\lam_2^U\otimes {1_{V_1}}^U$.

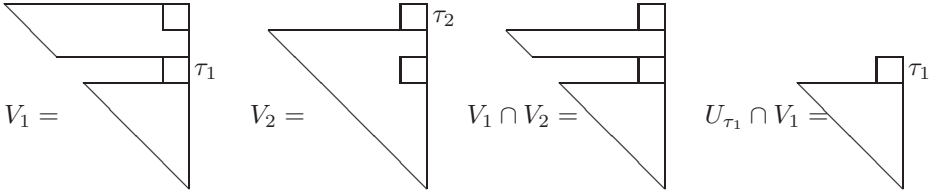
\begin{center}  

\begin{picture}(340,80)

\put(0,25){$V_1=$}
\put(0,70){\line(1,0){70}}  
\put(70,70){\line(0,-1){70}}  
\put(70,0){\line(-1,1){40}}   
\put(0,70){\line(1,-1){20}}  
\put(70,50){\line(-1,0){50}}
\put(70,40){\line(-1,0){40}}
\put(60,50){\line(0,-1){10}}
\put(60,70){\line(0,-1){10}}
\put(60,60){\line(1,0){10}}
\put(72,43){$\tau_1$}

\put(93,25){$V_2=$}
\put(150,70){\line(1,0){10}}  
\put(160,60){\line(-1,0){60}}
\put(150,70){\line(0,-1){10}}
\put(162,63){$\tau_2$}
\put(160,70){\line(0,-1){70}}  
\put(160,0){\line(-1,1){60}}   
\put(160,50){\line(-1,0){10}} \put(160,40){\line(-1,0){10}} \put(150,50){\line(0,-1){10}}

\put(175,25){$V_1\cap V_2=$}
\put(240,70){\line(1,0){10}}  
\put(250,60){\line(-1,0){60}}
\put(240,70){\line(0,-1){10}}
\put(250,70){\line(0,-1){70}}  
\put(250,0){\line(-1,1){40}}   
\put(190,60){\line(1,-1){10}}  
\put(250,50){\line(-1,0){50}} \put(250,40){\line(-1,0){40}} \put(240,50){\line(0,-1){10}}

\put(265,25){$U_{\tau_1}\cap V_1=$}
\put(340,50){\line(0,-1){50}}  
\put(340,0){\line(-1,1){40}}   
\put(340,50){\line(-1,0){10}} \put(340,40){\line(-1,0){40}} \put(330,50){\line(0,-1){10}}
\put(342,43){$\tau_1$}

\end{picture}
\end{center}

Here, ${1_{V_1}}^U=({1_{U_{\tau_1}\cap V_1}}^{U_{\tau_1}})_U$. Let
$M_0{:=}U_{\tau_1}\cap V_1$, and $M_t{:=}M_{t-1}X_{\alpha_{i,j-t}}$ for all
$t\in[1,j-i]$. Notice that $M_{j-i}=U_{\tau_1}$ and
${1_{U_{\tau_1}\cap V_1}}^{U_{\tau_1}}={1_{M_0}}^{U_{\tau_1}}$.

By the transitivity of inductions, ${1_{M_0}}^{U_{\tau_1}}$ is
decomposed via a series of inductions along the arm of $\tau_1$
namely from $M_0$ to $M_1$, $M_2$,..., $M_{j-i}$. Since
$M_1=M_0\rtimes X_{\alpha_{i,j-1}}$,
${1_{M_0}}^{M_1}=1_{M_1}+\sum_{\rho\in \Irr(M_1/M_0)^*}\rho$. For each
$\rho\in \Irr(M_1/M_0)^*$, by Lemma \ref{lem:midafi},
$\rho^{U_{\tau_1}}\in \Irr(U_{\tau_1})$ is an elementary character associated to
$\alpha_{i,j-1}$. Therefore,
$({1_{M_0}}^{M_1})^{U_{\tau_1}}={1_{M_1}}^{U_{\tau_1}}+\sum_{\rho\in
\Irr(M_1/M_0)^*}\rho^{U_{\tau_1}}$. The next induction step is from $M_1$ to $M_2$, i.e.
${1_{M_1}}^{U_{\tau_1}}=({1_{M_1}}^{M_2})^{U_{\tau_1}}$.

Similarly, $M_t{=}M_{t-1}{\rtimes} X_{\alpha_{i,j-t}}$, ${1_{M_{t-1}}}^{M_t}{=}1_{M_t}{+}\sum_{\rho\in
\Irr(M_t/M_{t-1})^*}\rho$, $\rho^{U_{\tau_1}}$ is an elementary
character associated to $\alpha_{i,j-t}$ for all $t\in[2,j-i]$.  So
${1_{M_0}}^{U_{\tau_1}}=1_{U_{\tau_1}}+\sum^{1\leq t\leq
j-i}_{\rho\in \Irr(M_t/M_{t-1})^*}\rho^{U_{\tau_1}}$.

Therefore, $({1_{M_0}}^{U_{\tau_1}})_U=1_U+ \sum^{1\leq t\leq
j-i}_{\rho\in \Irr(M_t/M_{t-1})^*}(\rho^{U_{\tau_1}})_U$. We have $H_{\alpha_{l,k}}\cap H_{\alpha_{i,j-t}}=\{1\}$ for all $t\in[1,j-i]$, by
(i) each component
${\lam_2}^U\otimes (\rho^{U_{\tau_2}})_U\in \Irr(U)$ and it appears in
$\lam_1^U\otimes\lam_2^U$ with multiplicity 1, for all above
$\rho$'s.

(iv)  Suppose $\alpha_{i,j}=\alpha_{l,k}$ and ${\lam_1}^U\neq
\overline{\lam_2}^U$. Hence, $\lam_1,\lam_2\in
\Irr(V_1/[V_1,V_1])^*$, and $\lam_1\neq \overline{\lam_2}$ where
$V_1=V_{\alpha_{i,j}}$. We define $V_2:=\la X_\alpha:\alpha\in \Sig^+-
leg(\alpha_{i,j}) \ra$. With the same argument as in Lemma
\ref{lem:midafi} on $leg(\alpha_{i,j})$, instead of $arm(\alpha_{i,j})$, there exists
$\xi\in \Irr(V_2/[V_2,V_2])^*$ such that
$\xi|_{X_{\alpha_{i,j}}}=\lam_2|_{X_{\alpha_{i,j}}}$,
$\xi|_{X_\alpha}=1_{X_\alpha}$ for all the others $X_\alpha\subset V_2$, and
$\xi^U={\lam_2}^U$ is an elementary character at $\alpha_{i,j}$.

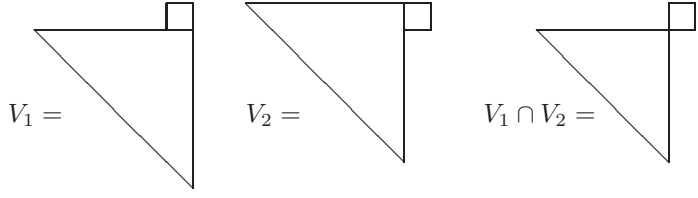
\begin{center} 
\begin{picture}(260,80)

\put(0,25){$V_1=$}
\put(60,70){\line(1,0){10}}  
\put(70,60){\line(-1,0){60}}
\put(60,70){\line(0,-1){10}}
\put(70,70){\line(0,-1){70}}  
\put(70,0){\line(-1,1){60}}   

\put(90,25){$V_2=$}
\put(90,70){\line(1,0){70}}  
\put(160,60){\line(-1,0){10}}
\put(160,70){\line(0,-1){10}}
\put(150,70){\line(0,-1){60}}  
\put(150,10){\line(-1,1){60}}   

\put(180,25){$V_1\cap V_2=$}
\put(250,70){\line(1,0){10}}  
\put(260,60){\line(-1,0){60}}
\put(260,70){\line(0,-1){10}}
\put(250,70){\line(0,-1){60}}  
\put(250,10){\line(-1,1){50}}   

\end{picture}
\end{center}

Let $V_{12}{:=}V_1\cap V_2$ and $\lam_{12}{:=}\lam_1\otimes \lam_2\in \Irr(V_1/[V_1,V_1])^*$. By Mackey formula for $U=V_2 V_1$,
\[\lam_1^U\otimes
{\xi}^U=({\lam_1}\otimes \xi^U|_{V_1})^U=(\lam_1\otimes {\xi|_{V_{12}}}^{V_1})^U=((\lam_1|_{V_{12}}\otimes
\xi|_{V_{12}})^{V_1})^U.\]
Since
$\xi|_{V_{12}}=\lam_2|_{V_{12}}$, $\lam_1^U\otimes
\xi^U=({\lam_{12}|_{V_{12}}}^{V_1})^U=(\lam_{12}\otimes
{1_{V_{12}}}^{V_1})^U$. Since $X_{\alpha_{i,j}}\subset
\ker({1_{V_{12}}}^{V_1})$, ${1_{V_{12}}}^{V_1}=(({1_{V_{12}}}^{V_1})_U)|_{ V_1}$. Therefore,
$\lam_1^U\otimes \lam_2^U=\lam_1^U\otimes {\xi}^U=\lam_{12}^U\otimes ({1_{V_{12}}}^{V_1})_U$. Here, $\lam_{12}^U$ is an elementary character associated to
$\alpha_{i,j}$ since $(\lam_1\otimes \lam_2)|_{X_{\alpha_{i,j}}}\neq
1_{X_{\alpha_{i,j}}}$.

We decompose the permutation character ${1_{V_{12}}}^{V_1}$
exactly with the same method as in (iii) into $1_{V_1}$ and a sum of
${\lam_\alpha}^U$, for all $\alpha\in arm(\alpha_{i,j})$ and all $\lam_\alpha\in \Irr(V_\alpha/[V_\alpha,V_\alpha])^*$.

So $\lam_1^U\otimes \lam_2^U=\lam_{12}^U+\lam_{12}^U\otimes\sum^{\alpha\in arm(\alpha_{i,j})}_{\lam_\alpha\in \Irr(V_\alpha/[V_\alpha,V_\alpha])^*}\lam_\alpha^U$. At the
simple root $\alpha_i\in arm(\alpha_{i,j})$, as observed above, $\lam_{12}^U\otimes\rho=\lam_{12}^U$ for all linear
elementary $\rho$ associated to $\alpha_i$. Hence, we obtain $q$ copies of $\lam_{12}^U$. At others
 $\alpha_{i,t}$ where $i<t<j$, by (iii) with $\alpha_{i,t}\in arm(\alpha_{i,j})$, each tensor product $\lam_{12}^U\otimes \lam_{\alpha_{i,t}}^U$ decomposes into  $\lam_{12}^U$ and a sum of $\sum^{\alpha\in leg(\alpha_{i,t})}_{\lam_\alpha\in \Irr(V_\alpha/[V_\alpha,V_\alpha])^*}\lam_{12}^U\otimes \lam_\alpha^U$. This gives the others $(j-i-1)(q-1)$  $\lam_{12}^U$. All the other tensor products are irreducible by (i).

(v) Suppose $\alpha_{i,j}=\alpha_{l,k}$ and
$\lam_1^U=\overline{\lam_2^U}$, which implies
$\lam_1=\overline{\lam_2}$. We use the same setup as in (iv), hence
$\lam_2|_{V_{12}}=\xi|_{V_{12}}=\overline{\lam_1}|_{V_{12}}$, and
\[\lam_1^U\otimes \xi^U=({\lam_1}|_{V_{12}}\otimes \xi|_{V_{12}})^U={1_{V_{12}}}^U=({1_{V_{12}}}^{V_1})^U.\]
We use exactly the same method as in (iii) to decompose the permutation
character ${1_{V_{12}}}^{V_1}$ into $1_{V_1}$ and a sum of all
elementary characters at all $\alpha\in arm(\alpha_{i,j})$, each
constituent appears with multiplicity 1. Call this set of
constituents $A$.

Since each character $\xi\in A$ is inflated to $U$,
$\xi^U=\xi_U\otimes {1_{V_2}}^U$. Again we decompose ${1_{V_2}}^U$
into $1_U$ and a sum of all elementary characters at all $\alpha\in
leg(\alpha_{i,j})$, each constituent appears with multiplicity 1. Call
this set of constituents $B$.

Therefore, $\lam_1^U\otimes \overline{\lam_1}^U=({1_{V_{12}}}^{V_1})^U=\sum_{\xi_1\in A}{\xi_1}^U=\sum_{\xi_1\in A}\sum_{\xi_2\in B}\ (\xi_1)_U\otimes \xi_2$.
\end{proof}

\medskip
Notice that in the proof of Lemma \ref{lem:prod2midafis} (ii) because $V\ntrianglelefteq U$, we use a sequence of inductions from $M_t$ to $M_{t+1}\subset N_U(M_t)$ to control distinct induced characters. This technique will be reused in our latter proofs. Let $1_U$  be the supercharacter $\xi_{D,\phi}$ when $D$ is empty. Lemma \ref{lem:prod2midafis} (iii), (iv),
(v) give the proof of the following.

\begin{corollary} \label{cor:tensor=sum_of_tensors}

Let $S\subset \Sig^+$ be nonempty. The tensor product of a set of elementary characters at $\alpha\in S$ decomposes into a sum of supercharacters.

\end{corollary}

By the proof in Lemma \ref{lem:prod2midafis}, it is easy to
count the multiplicity of each supercharacter constituent of the tensor product in Corollary
\ref{cor:tensor=sum_of_tensors}. Now, we prove Theorem \ref{thm:Andre}.

\medskip
\begin{proof}[Proof of Theorem \ref{thm:Andre}] We show that the regular character ${1_{\{1\}}}^U$  decomposes into $1_U$ and a
sum of supercharacters by induction on
$n$, where $U=U_n(q)$, and that there is no $\chi\in \Irr(U)^*$ such that $\chi$ is
an irreducible constituent of two distinct supercharacters $\xi_{D,\phi}$,
$\xi_{D',\phi'}$.

Let $M{:=}X_{\alpha_{1,n-1}}L_{\alpha_{1,n-1}}$, $S{:=}\{\alpha_{1,n-1}\}\cup
leg(\alpha_{1,n-1})$. Since $M{=}\la
X_\alpha:\alpha\in S\ra$ is abelian,
${1_{\{1\}}}^M=\sum_{\rho\in \Irr(M)}\rho
=1_M+\sum^{\alpha\in S}_{\lam\in \Irr(V_\alpha/[V_\alpha,V_\alpha])^*}\lam|_M\otimes ({1_{\{1\}}}^{L_\alpha})_M$.
Here ${1_{\{1\}}}^{L_\alpha}$ is the regular character of $L_\alpha$, by
Corollary \ref{cor:midafi_to_leg}, $\lam|_M\otimes
({1_{\{1\}}}^{L_\alpha})_M=(\lam^U|_{X_\alpha L_\alpha})_M=\lam^U|_M$ for all
$\lam\in \Irr(V_\alpha/[V_\alpha,V_\alpha])^*$ and $\alpha\in S$. Therefore,
${1_{\{1\}}}^M=1_M+\sum^{\alpha\in S}_{\lam\in \Irr(V_\alpha/[V_\alpha,V_\alpha])^*}\lam^U|_M$.
By the transitivity of inductions,
\[\begin{array}{lll}
{1_{\{1\}}}^U&=&({1_{\{1\}}}^M)^U\\
&=&(1_M+\sum^{\alpha\in S}_{\lam\in \Irr(V_\alpha/[V_\alpha,V_\alpha])^*}\lam^U|_M)^U\\
&=&{1_M}^U+\sum^{\alpha\in S}_{\lam\in \Irr(V_\alpha/[V_\alpha,V_\alpha])^*}\lam^U\otimes {1_M}^U.
\end{array}\]

Since $U=M\rtimes U_{\alpha_{1,n-2}}$, ${1_M}^U$ is the regular character
of $U_{\alpha_{1,n-2}}$ inflated to $U$. By induction on $n$, ${1_M}^U$ equals $1_U$ and a sum of
supercharacters of $U_{\alpha_{1,n-2}}$. By
Corollary \ref{cor:tensor=sum_of_tensors}, ${1_{\{1\}}}^U$ equals $1_U$ and a
sum of supercharacters.

About multiplicities, the decomposition of ${1_{\{1\}}}^U$ into $1_U$ and a sum of supercharacters shows that $1_U$ is not a constituent of any supercharacters.

\medskip
(Uniqueness) Suppose the set of constituents of supercharacters does not partition $\Irr(U)^*$. There is $\chi\in \Irr(U)^*$ such that
$\chi$ is a constituent of distinct supercharacters $\xi_{D,\phi}$ and $\xi_{D',\phi'}$. It is clear that $|D|>1$, and $|D|$ can be chosen to be the smallest. Since $\xi_{D,\phi_D}$ and $\xi_{D',\phi_{D'}}$ have at least a common constituent,  $(\xi_{D,\phi},\xi_{D',\phi'})=(\xi_{D,\phi}\otimes\overline{\xi_{D',\phi'}},1_U)\neq 0$.

We start with a highest root in $D\cup D'$, call it $\alpha$. Suppose $\alpha\in D$ and we shall show that $\alpha\in D'$. We have $X_\alpha\subset
Z(\xi_{D,\phi}\otimes\overline{\xi_{D',\phi'}})$. If $\alpha\notin D'$, then $X_\alpha\subset Z(\psi)$ and $X_{\alpha}\nsubseteq \ker(\psi)$ for all
constituents $\psi$ of
$\xi_{D,\phi}\otimes\overline{\xi_{D',\phi'}}$, which implies
$(\xi_{D,\phi}\otimes\overline{\xi_{D',\phi'}},1_U)=0$, contradictory to our assumption. So $\alpha\in D'$. The corresponding characters
$\lam_\alpha\in \phi$ and  ${\lam'}_\alpha\in \phi'$ must satisfy $(\lam_\alpha|_{X_\alpha}\otimes \overline{\lam'_\alpha}|_{X_\alpha},1_{X_\alpha})=1$, i.e. $\lam_\alpha=\lam'_\alpha$; otherwise, by Lemma
\ref{lem:prod2midafis} (iv), $X_\alpha\nsubseteq \ker(\psi)$ for all constituents of $\xi_{D,\phi}\otimes\overline{\xi_{D',\phi'}}$.

By Lemma \ref{lem:prod2midafis} (v), $\lam_\alpha^U\otimes \overline{\lam_\alpha}^U$ decomposes into $1_U$
and a sum of $\lam_{\beta_1}^U$, $\lam_{\beta_2}^U$, $\lam_{\beta_1}^U\otimes
\lam_{\beta_2}^U$ for all $\beta_1\in arm(\alpha)$, $\beta_2\in
leg(\alpha)$ and all $\lam_{\beta_k}\in \Irr(V_{\beta_k}/[V_{\beta_k},V_{\beta_k}])^*$ with $k=1,2$. Since $D'$ is basic, for each $\beta_k$, if ${\beta_k}\cup (D'- \{\alpha\})$ is not basic then there is some $\gamma\in (D'-\{\alpha\})$ such that either $\beta_k\in arm(\gamma)$ (or $leg(\gamma)$) or $\gamma\in arm(\beta_k)$ (or $leg(\beta_k)$).
By Corollary \ref{cor:tensor=sum_of_tensors} and Lemma \ref{lem:prod2midafis} (vi), (v), $\overline{\lam_\alpha}^U\otimes\xi_{D',\phi'}$ decomposes into a sum of supercharacters $\xi_{S,\phi_S}$. Corresponding to constituents of $\lam_\alpha^U\otimes \overline{\lam_\alpha}^U$, since $(\phi- \{\lam_\alpha\})\neq (\phi'-\{\lam_\alpha\})$ and others $S\neq (D- \{\alpha\})$ because $S$ contains at least a root in $arm(\alpha)$ or $leg(\alpha)$, none of constituents equal $\bigotimes_{\lam_\delta\in\phi-\{\lam_\alpha\}} \lam_\delta^U$.

Since $0\neq (\xi_{D,\phi},\xi_{D',\phi'})=(\bigotimes_{\lam_\delta\in \phi}{\lam_\delta}^U,\xi_{D',\phi'})=(\bigotimes_{\lam_\delta\in \phi-\{\lam_\alpha\}} \lam_\delta^U,\overline{\lam_\alpha}^U\otimes \xi_{D',\phi'})$, there exists some $\xi_{S_0,\phi_{S_0}}$ of $\overline{\lam_\alpha}^U\otimes\xi_{D',\phi'}$ such that
$(\bigotimes_{\lam_\delta\in \phi,\delta\neq \alpha} \lam_\delta^U,\xi_{S_0,\phi_{S_0}})\neq 0$, contradictory to the minimality of $|D|$ since $D- \{\alpha\}$ is basic.
\end{proof}

\smallskip
By the uniqueness in Theorem \ref{thm:Andre}, the constituent set of
all supercharacters $\xi_{D,\phi}$ is a partition of $\Irr(U)^*$.
We define that a basic set $D$ is {\it decomposable} if $D$ has
two nonempty disjoint subsets $A$, $B$ such that $A\cup B=D$, $\la
U_\alpha:\alpha\in A \ra \subset U_\beta$ and  $\la H_\alpha:\alpha\in B
\ra\subset R_\beta$ for some $\beta\in \Sig^+$; a basic set $D$
is {\it indecomposable} if it is not decomposable.

For example, in $U_7(q)$,
$D=\{\alpha_{1,2},\alpha_{3,4},\alpha_{2,5},\alpha_{4,6}\}$ is indecomposable,
and $D=\{\alpha_3,\alpha_4\}\cup \{\alpha_{2,5},\alpha_{1,6}\}$ is decomposable.
Denote by $\Irr(\xi_{D,\phi})$ the irreducible constituent set
of $\xi_{D,\phi}$. The following property is very
helpful to work with supercharacters.

\begin{lemma} \label{lem:D_decomposable}

Let $D$ be a decomposable basic set as $A\cup B$. Then for each $\phi\in E_D$, $\phi_A\in E_A$, $\phi_B\in E_B$ such that $\phi=\phi_A\cup\phi_B$, $\Irr(\xi_{D,\phi})=\{\chi_A\otimes\chi_B:\chi_A\in \Irr(\xi_{A,\phi_A}),\ \chi_B\in \Irr(\xi_{B,\phi_B})\}$, i.e. $\Irr(\xi_{D,\phi})\simeq \Irr(\xi_{A,\phi_A})\times \Irr(\xi_{B,\phi_B})$.

\end{lemma}

\begin{proof}
Suppose that $D=A\cup\!\!\!\!\cdot\ B$
where $\la U_\alpha:\alpha\in A \ra \subset U_\beta$ and  $\la
H_\alpha:\alpha\in B \ra\subset R_\beta$ for some $\beta\in \Sig^+$.
For $\phi\in E_D$, let $\phi_A=\{\lam_\alpha\in \phi:\alpha\in A\}\in
E_A$ and $\phi_B=\{\lam_\alpha\in \phi:\alpha\in B\}\in E_B$. By the
uniqueness in Theorem \ref{thm:Andre},  for all
$\chi_A\in \Irr(\xi_{A,\phi_A})$, $\chi_B\in \Irr(\xi_{B,\phi_B})$, we
show that $\chi_A\otimes\chi_B\in \Irr(U)$. It is clear that $\chi|_{U_\beta}\in \Irr(U_\beta)$ for all $\chi\in \Irr(\xi_{A,\phi_A})$.  Since $\la
H_\alpha:\alpha\in B \ra \subset R_\beta\unlhd U$,
by Clifford theory, it suffices to show that $\xi_{B,\phi_B}|_{R_\beta}$ has the same number
of constituents of $\xi_{B,\phi_B}$, i.e. for all $\chi\in
\Irr(\xi_{B,\phi_B})$, $\chi|_{R_\beta}\in \Irr(R_\beta)$, because the
number of irreducible constituents of $\xi_{B,\phi_B}|_{R_\beta}$ is
greater than or equal to the one of
$\xi_{B,\phi_B}$. Therefore,  we need to show that
$(\chi,\xi_{B,\phi_B})=(\chi|_{R_\beta},\xi_{B,\phi_B}|_{R_\beta})$
for all $\chi\in \Irr(\xi_{B,\phi_B})$. By Frobenius reciprocity,
\[\begin{array}{lll}
(\chi|_{R_\beta},\xi_{B,\phi_B}|_{R_\beta})&=&(\chi,(\xi_{B,\phi_B}|_{R_\beta})^U)\\
&=&(\chi,(\xi_{B,\phi_B}|_{R_\beta}\otimes 1_{R_\beta})^U)\\
&=&(\chi,\xi_{B,\phi_B}\otimes{1_{R_\beta}}^U).
\end{array}\]

Since $U=R_\beta\rtimes U_\beta$, ${1_{R_\beta}}^U=({1_{\{1\}}}^{U_\beta})_U$ is the regular
character of $U_\beta$ inflated to $U$. By the proof of Theorem \ref{thm:Andre},
${1_{R_\beta}}^U$ equals $1_U$ and a sum of all supercharacters of $U_{\beta}$. Since $B\cup
T$ is a basic set of $U$ for all basic sets  $T$ of $U_\beta$, by
Theorem \ref{thm:Andre}, the uniqueness of $\chi\in
\Irr(\xi_{B,\phi_B})$ forces
$(\chi,\xi_{B,\phi_B}\otimes{1_{R_\beta}}^U)=(\chi,\xi_{B,\phi_B})$.
\end{proof}


\section{Irreducible representations of $U$ of large degree} \label{section:three_largest}

It is known  in \cite{Is2,Hup} that the degrees of all irreducible characters of
$U$ are precisely $\{q^e|\ 0\leq e\leq \mu(n)\}$ where \[\mu(n)=\left\{
\begin{array}{ll} m(m-1)& \mbox{if } n=2m,  \\ m^2 & \mbox{if } n=2m+1. \end{array} \right.\]

For each $0\leq e \leq \mu(n)$,  denote by $N_{n,e}(q)$ the number of irreducible characters of $U_n(q)$
having degree $q^e$.
We need some more preparation before presenting the construction of large degree irreducible representations.

\begin{lemma} \label{p1}

For $1\leq k <n-1$, let $S_k:=\Sig^+-
(\{\alpha_{1,k},\alpha_{k+1,n-1}\}\cup leg(\alpha_{1,k})\cup
arm(\alpha_{k+1,n-1}))$. The following are true:
\begin{itemize}
\item[(i)] $S_k$ is closed and $T_k:=\la X_\alpha:\alpha\in S_k\ra$
is isomorphic to $U_{n-1}(q)$.
\item[(ii)] $[T_k,X_{\alpha_{1,k}}]=[T_k,X_{\alpha_{k+1,n-1}}]=\{1\}$.
\end{itemize}

\end{lemma}

\begin{proof}
(i) Let $\{\beta_i:1\leq i\leq n-2\}$ be the positive
simple root set of $\Sigma_{n-2}$. Let $\varphi$ be the function
mapping $\alpha_i$ to $\beta_i$ for $i\in[1,k-1]$,  $\alpha_{k,k+1}$ to
$\beta_k$, and $\alpha_{i+1}$ to $\beta_{i}$ for $i\in[k+1,n-2]$. The
claim is clear by extending $\varphi$ linearly to all $S_k$ and mapping $x_\alpha(c)$ to $x_{\varphi(\alpha)}(c)$ for
all $\alpha\in S_k$ and $c\in \F_q$.

(ii) Since $[X_\alpha,X_{\alpha_{1,k}}]=[X_\alpha,X_{\alpha_{k+1,n-1}}]=\{1\}$ for all $\alpha\in S_k$, the statement is clear.
\end{proof}

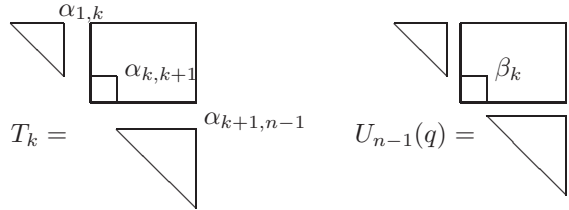
\begin{center}  

\begin{picture}(210,80)

\put(0,25){$T_k=$}
\put(0,70){\line(1,0){20}}   \put(30,70){\line(1,0){40}}  
\put(70,70){\line(0,-1){30}} \put(70,30){\line(0,-1){30}} 
\put(70,30){\line(-1,0){30}} \put(70,40){\line(-1,0){40}}
\put(20,70){\line(0,-1){20}} \put(30,70){\line(0,-1){30}}
\put(70,0){\line(-1,1){30}}  \put(0,70){\line(1,-1){20}}   

\put(18,73){$\alpha_{1,k}$}
\put(73,32){$\alpha_{k+1,n-1}$}
\put(43,50){$\alpha_{k,k+1}$}
\put(30,50){\line(1,0){10}} \put(40,50){\line(0,-1){10}}

\put(130,25){$U_{n-1}(q)=$}
\put(145,70){\line(1,0){20}}   \put(170,70){\line(1,0){40}}  
\put(210,70){\line(0,-1){30}} \put(210,35){\line(0,-1){30}} 
\put(210,35){\line(-1,0){30}} \put(210,40){\line(-1,0){40}}
\put(165,70){\line(0,-1){20}} \put(170,70){\line(0,-1){30}}
\put(210,5){\line(-1,1){30}}  \put(145,70){\line(1,-1){20}}   
\put(170,50){\line(1,0){10}} \put(180,50){\line(0,-1){10}}
\put(183,50){$\beta_k$}

\end{picture}
\end{center}

It is clear that $T_k\cap U_{\alpha_{1,n-t}}\simeq U_{n-t}(q)$ for all $t\in[1,n-k-1]$. To facilitate the use of Lemma \ref{p1}, we make the following.

\begin{definition}
We define $T_{k,t}:= T_k\cap U_{\alpha_{1,n-t}}$ where $T_k$ is in Lemma \ref{p1}.
\end{definition}

For each almost faithful $\chi\in \Irr(U)$, by Theorem \ref{thm:unique_factor}, there exist uniquely $\psi\in \Irr(U_{\alpha_{2,n-2}})$ and $\lam\in \Irr(V_{\alpha_{1,n-1}}/[V_{\alpha_{1,n-1}},V_{\alpha_{1,n-1}}])^*$ such that $\chi=(\psi_{V_{\alpha_{1,n-1}}}\otimes \lam)^U$. The next lemma generalizes this construction for all irreducible characters.

\begin{lemma} \label{p2}

Set $\alpha=\alpha_{k+1,n-1}$ for some $k\in[1,n-2]$. Suppose that there is $\chi\in \Irr(U)$ such that $X_{\alpha}\subset Z(\chi)$ and $X_\alpha\nsubseteq \ker(\chi)$.  Then there exist $\xi\in \Irr(V_\alpha\cap U_{\alpha_{1,n-2}})$ and $\lambda\in \Irr(V_\alpha/[V_\alpha,V_\alpha])^*$ such that $\chi=(\xi_{V_\alpha}\otimes \lambda)^U$.

\end{lemma}

\begin{proof}
Set $H=V_\alpha\cap U_{\alpha_{1,n-2}}$. Using the same method as in Lemma \ref{lem:prod2midafis} (ii) by applying a sequence of inductions along the arm of $\alpha$, $(\xi_{V_\alpha}\otimes \lambda)^U\in \Irr(U)$ for all $\xi\in \Irr(H)$ and $\lambda\in \Irr(V_\alpha/[V_\alpha,V_\alpha])^*$. Since $X_{\alpha}\subset Z(\chi)$ and $X_\alpha\nsubseteq \ker(\chi)$, there exists $\lambda\in \Irr(V_\alpha/[V_\alpha,V_\alpha])^*$ such that $\chi(x)=\lambda(x)\chi(1)$ for all $x\in X_\alpha$. Hence, we shall show that $\chi\in \Irr(U,\xi_{V_\alpha}\otimes \lambda)$ for some $\xi\in \Irr(H)$.

Let $\phi$ be a representation affording $\chi$. For each $t\in[1,k]$, all $x\in X_{\alpha}$ and $y\in X_{\alpha_{t,k}}$,
$\phi([y,x]){=}[\phi(y),\phi(x)]{=}[\phi(y),\lambda(x)\phi(1)]{=}\phi(1)$. Since $X_{\alpha_{t,n-1}}{=}[X_{\alpha_{t,k}},X_\alpha]$, $X_{\alpha_{t,n-1}}\subset \ker(\chi)$. We proceed $U$ by modulo $\prod_{t=1}^kX_{\alpha_{t,n-1}}$. So $L_\alpha X_\alpha \unlhd U$ where $L_\alpha=\la X_\beta: \ \beta \in leg(\alpha)\ra$.

Since the hook group $H_\alpha$ is special,  $\chi|_{H_\alpha}=a\cdot \lambda^U|_{H_\alpha}$ for some $a\in \Z^+$. By Corollary \ref{cor:midafi_to_leg}, $\chi|_{L_\alpha}=a\cdot {1_{\{1\}}}^{L_\alpha}$. Thus, $\chi\in \Irr(U,1_{L_\alpha})$ because $(\chi|_{L_\alpha},1_{L_\alpha})=(a\cdot {1_{\{1\}}}^{L_\alpha},1_{L_\alpha})=a>0$.

 Suppose that $\chi|_{V_\alpha}=\sum_{\xi\in S}a_\xi\cdot \xi$ for some $S\subset \Irr(V_\alpha,\lambda|_{X_\alpha})$ and $a_\xi\in \Z^+$. Since  $\chi\in \Irr(U,1_{L_\alpha})$, there exists $\xi_0\in S$ such that $\xi_0\in \Irr(V_\alpha,1_{L_\alpha})$. So $\xi_0\in \Irr(V_\alpha,\lambda|_{L_\alpha X_\alpha})$.  Since $L_\alpha X_\alpha \unlhd V_\alpha$ and $V_\alpha=L_\alpha X_\alpha\rtimes H$, ${1_{L_\alpha X_\alpha}}^{V_\alpha}= ({1_{\{1\}}}^H)_{V_\alpha}$, the regular character of $H$ inflated to $V_\alpha$. Thus, ${\lambda|_{L_\alpha X_\alpha}}^{V_\alpha}=\lambda\otimes ({1_{\{1\}}}^H)_{V_\alpha}$. Since $\xi_0\in \Irr(V_\alpha,\lambda|_{L_\alpha X_\alpha})$, $\xi_0|_H\in \Irr(H)$ and $\xi_0=\lambda\otimes(\xi_0|_H)_{V_\alpha}$. So $\chi\in \Irr(U,\xi_0)=\{{\xi_0}^U\}$.
 \end{proof}

\smallskip
Actually the uniqueness of $\xi$ and $\lam$ in Lemma \ref{p2} can be proved. Now we present the largest degree irreducible representations of $U$ which has been stated in Theorem \ref{intro:main-thm}.

\subsection{The largest degree irreducible representations of $U$}
 Lehrer \cite{Leh} shows that the degree of a largest degree irreducible character is $q^{(n-2)+(n-4)+...}=q^{\mu(n)}$, and that the number of irreducible characters of degree $q^{\mu(n)}$ is at least $(q -1)^{[n/2]}$. (The square brackets here denote the `greatest
integer' function.) By using Lehrer's regular orbit method \cite{Leh}, $\mu(n)>\mu(n-1)$ for all $n\geq 3$.

\begin{slemma} \label{slem:largest}

For $n\geq 3$, each largest degree irreducible character of $U$ factors uniquely  into a tensor product of an almost faithful elementary character with a largest degree irreducible character of $U_{\alpha_{2,n-2}}$. Therefore, \[N_{n,\mu(n)}(q)=(q-1)N_{n-2,\mu(n-2)}(q).\]

\end{slemma}

\begin{proof}
By Theorem \ref{thm:unique_factor}, a product of an almost faithful elementary character and a largest degree irreducible character of $U_{\alpha_{2,n-2}}$ is almost faithful and has degree equal to $q^{n-2+\mu(n-2)}$. Since the described ones have the largest degree in the set of all almost faithful irreducibles, it is enough to show that not almost faithful $\chi\in \Irr(U)$ has degree less than $q^{n-2+\mu(n-2)}$. If $X_{\alpha_{k,n-1}}\subset \ker(\chi)$ for all $1\leq k\leq n-1$, $\chi$ can be considered as a character of $U_{\alpha_{1,n-2}}$ and hence, $\chi(1)\leq q^{\mu(n-1)}<q^{\mu(n)}$ since $n\geq 3$. Therefore, we suppose that there is some $1\leq k<n-1$ such that $X_{\alpha_{k+1,n-1}}\subset Z(\chi)$ and $X_{\alpha_{k+1,n-1}}\nsubseteq \ker(\chi)$. Set $\alpha=\alpha_{k+1,n-1}$.

By Lemma \ref{p2}, $\chi=(\xi_{V_\alpha}\otimes \lambda)^U$ for some $\xi\in \Irr(V_\alpha\cap U_{\alpha_{1,n-2}})$ and $\lambda\in \Irr(V_\alpha/[V_\alpha,V_\alpha])^*$.  Since $\chi(1)=\xi(1)\lambda^U(1)=q^{n-k-2}\xi(1)$,  we shall show that $\xi(1)<q^{\mu(n-2)+k}$.

Let $T{:=}\la X_\beta\subset V_\alpha\cap U_{\alpha_{1,n-2}}:\ \beta\notin leg(\alpha_{1,k})\ra = T_{k,2}{\times} X_{\alpha_{1,k}}\simeq U_{n-2}(q){\times} X_{\alpha_{1,k}}$ by Lemma \ref{p1}. So a largest degree irreducible character of $T$ has degree $q^{\mu(n-2)}$. Since $|V_\alpha:T|=q^{|leg(\alpha_{1,k})|}=q^{k-1}$ and each irreducible character of $V_\alpha$ is a constituent of an induced character of $T$, $\xi(1)\leq q^{\mu(n-2)+k-1}<q^{\mu(n-2)+k}$.
\end{proof}

\smallskip
From the proof of Lemma \ref{slem:largest}, the recursion formula $\mu(n)=n-2+\mu(n-2)$ is directly obtained. We can obtain back the formula of $\mu(n)$ as stated above with the initial values $\mu(2)=0$ and $\mu(3)=1$. The following gives the connection between $\mu(n)$ and $\mu(n-1)$.

\begin{slemma}  \label{slem:mu_n-1}

$\mu(n)-\mu(n-1)=[(n-1)/2]$ for $n\geq 2$.

\end{slemma}
\begin{proof}
This is Lemma 7.2 (a) in \cite{Is3}.
\end{proof}

\smallskip
Our proof shows that the largest degree irreducible characters are supercharacters $\xi_{D,\phi}$ where $D=\{\alpha_{k,n-k}: 1\leq k\leq m\}$ if $n=2m+1$, and $D=\{\alpha_{k,n-k}: 1\leq k\leq m\}$ or $D=\{\alpha_{k,n-k}: 1\leq k< m\}$ if $n=2m$, because at $\alpha_{m,2m-m}=\alpha_m$, all elementary characters are linear. Notice that
Isaacs \cite{Is3} computes the number of regular orbits and obtains the recursion formula $N_{n,\mu(n)}(q)$ in Lemma \ref{slem:largest} in his proof. Moreover, he shows
\[N_{n,\mu(n)}(q)=\left\{\begin{array}{ll} (q-1)^m &\mbox{ if } n=2m+1,  \\q(q-1)^{m-1} &\mbox{ if } n=2m. \end{array} \right.\]
In another way, we could use this result to prove Lemma \ref{slem:largest} by applying Theorem \ref{thm:Andre} and counting the supercharacters described above.

\subsection{The second largest degree irreducible representations of $U$}

There are two cases as follows.
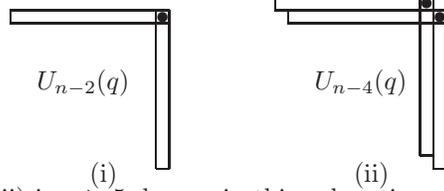
\begin{center} 

\begin{picture}(165,65)

\put(0,60){\line(1,0){60}} \put(0,55){\line(1,0){60}} \put(0,60){\line(0,-1){5}}
\put(60,60){\line(0,-1){60}}  \put(55,60){\line(0,-1){60}} \put(55,0){\line(1,0){5}}
\put(55,55){$\bullet$}
\put(10,30){$U_{n-2}(q)$}
\put(30,-5){(i)}

\put(100,65){\line(1,0){60}}  \put(100,60){\line(1,0){60}} \put(100,65){\line(0,-1){5}}
\put(160,65){\line(0,-1){60}}  \put(155,65){\line(0,-1){60}} \put(155,5){\line(1,0){5}}
\put(155,60){$\bullet$}
\put(105,60){\line(1,0){60}}  \put(105,55){\line(1,0){60}} \put(105,60){\line(0,-1){5}}
\put(165,60){\line(0,-1){60}}  \put(160,60){\line(0,-1){60}} \put(160,0){\line(1,0){5}}
\put(160,55){$\bullet$}
\put(115,30){$U_{n-4}(q)$}
\put(130,-5){(ii)}
\end{picture}

\end{center}

The condition of case (ii) is $n\geq 5$, hence, in this subsection we suppose $n\geq 5$.

In (i), by Theorem
\ref{thm:unique_factor}, each tensor product of an almost faithful elementary character with a second
largest degree irreducible character of $U_{\alpha_{2,n-2}}\simeq
U_{n-2}(q)$ is a second largest degree  almost faithful irreducible
character of $U$. There are $(q-1)N_{n-2,\mu(n-2)-1}(q)$ second largest degree irreducible characters
in this case.

In case (ii), each tensor product of two elementary characters associated to
$\alpha_{1,n-2}$ and $\alpha_{2,n-1}$ respectively, by Lemma
\ref{lem:prod2midafis} (ii), decomposes into $q$ distinct
irreducible constituents of degree $q^{n-3+n-3-1}=q^{2n-7}$.
By Lemma \ref{lem:D_decomposable}, a tensor product of
each constituent with a largest degree irreducible character of
$U_{\alpha_{3,n-3}}\simeq U_{n-4}(q)$ is a second largest irreducible character
of $U$. Since
$X_{\alpha_{1,n-1}}$ is in their kernels, these characters are not almost faithful. So we obtain
$q(q-1)^2N_{n-4,\mu(n-4)}(q)$ second largest degree irreducible characters in this case.

\begin{slemma}\label{slem:second_largest}

For $n\geq 5$, each second largest degree irreducible representation
of $U$ is of the form (i) or (ii). Therefore,
\[N_{n,\mu(n)-1}(q)=(q-1)N_{n-2,\mu(n-2)-1}(q)+q(q-1)^2N_{n-4,\mu(n-4)}(q).\]

\end{slemma}

\begin{proof}
By Theorem \ref{thm:unique_factor} and Lemma \ref{lem:D_decomposable}, it suffices to show that $\chi(1)<q^{\mu(n)-1}=q^{n-3+\mu(n-2)}$ for all $\chi\in \Irr(U)$ such that $X_{\alpha_{1,n-2}}X_{\alpha_{1,n-1}}\subset \ker(\chi)$.  Notice that $X_{\alpha_{1,n-2}}X_{\alpha_{1,n-1}}$ is mapped to  $X_{\alpha_{2,n-1}}X_{\alpha_{1,n-1}}$ by the graph automorphism of $\GL_n(q)$, i.e. the map sending $g$ to ${w_0}^{-1}g^Tw_0$ where $w_0$ is the longest element in the Weyl group $S_n$ and $T$ means the transpose. We proceed $U$ by modulo $X_{\alpha_{1,n-2}}X_{\alpha_{1,n-1}}$.
If $X_{\alpha_{k,n-1}}\subset \ker(\chi)$ for all $i\in[1,n-1]$, $\chi$ can be considered as a character of $U_{\alpha_{1,n-2}}$, hence, $\chi(1)\leq q^{\mu(n-1)}<q^{\mu(n)-1}$ by Lemma \ref{slem:mu_n-1} with $n\geq 5$. Therefore, we suppose that there is $k\in[1,n-2]$ such that $X_{\alpha_{k+1,n-1}}\subset Z(\chi)$ and $X_{\alpha_{k+1,n-1}}\nsubseteq \ker(\chi)$.  Set $\alpha:=\alpha_{k+1,n-1}$.

Now we apply the same technique as in Lemma \ref{slem:largest}. By Lemma \ref{p2}, $\chi=(\xi_{V_\alpha}\otimes \lambda)^U$ for some $\xi\in \Irr(V_\alpha\cap U_{\alpha_{1,n-2}})$ and $\lambda\in \Irr(V_\alpha/[V_\alpha,V_\alpha])^*$. Since $\chi(1)=\xi(1)\lambda^U(1)=q^{n-k-2}\xi(1)$, we shall show that $\xi(1)<q^{\mu(n-2)+k-1}$.

Let $T:=\la X_\beta\subset V_\alpha \cap U_{\alpha_{1,n-2}}: \beta\notin leg(\alpha_{1,k})\ra = T_{k,2}\times X_{\alpha_{1,k}}\simeq U_{n-2}(q)\times X_{\alpha_{1,k}}$ by Lemma \ref{p1}. Since we work with $U/X_{\alpha_{1,n-2}}X_{\alpha_{1,n-1}}$,  $T_{k,2}$ does not have any almost faithful characters. So a largest degree irreducible character of $T$ has degree at most $q^{\mu(n-2)-1}$. Since $|V_\alpha:T|=q^{|leg(\alpha_{1,k})|}=q^{k-1}$ and each irreducible character of $V_\alpha$ is a constituent of some induced character of $T$, $\xi(1)\leq q^{\mu(n-2)-1+k-1}<q^{\mu(n-2)+k-1}$.
\end{proof}

\smallskip
As above, we remark that by using combinatorial methods, Isaacs \cite{Is3} obtains the recursion formula in Lemma \ref{slem:second_largest} in his proof, and then, calculates
\[N_{n,\mu(n)-1}(q)=\left\{\begin{array}{cl} q(q-1)^{m-1}(m(q-1)+1)&\mbox{ if } n=2m+1, \\ q(q-1)^{m-1}((m-1)q+1)&\mbox{ if } n=2m. \end{array}\right.\]
Again we could use this result to obtain all second largest degree irreducible representations in Lemma \ref{slem:second_largest} by counting irreducibles described in cases (i) and (ii).

\subsection{The third largest degree irreducible representations of $U$}

There are five cases as follows.
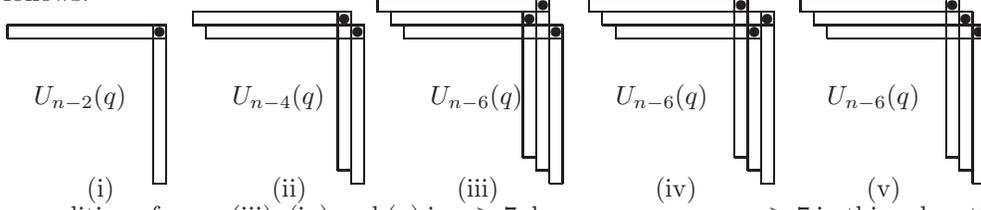
\begin{center} 
\begin{picture}(370,65)

\put(0,60){\line(1,0){60}} \put(0,55){\line(1,0){60}} \put(0,60){\line(0,-1){5}}
\put(60,60){\line(0,-1){60}}  \put(55,60){\line(0,-1){60}} \put(55,0){\line(1,0){5}}
\put(10,30){$U_{n-2}(q)$}
\put(55,55){$\bullet$}
\put(30,-5){(i)}

\put(70,65){\line(1,0){60}}   \put(70,60){\line(1,0){60}}    \put(70,65){\line(0,-1){5}}
\put(130,65){\line(0,-1){60}}  \put(125,65){\line(0,-1){60}} \put(125,5){\line(1,0){5}}
\put(125,60){$\bullet$}
\put(75,60){\line(1,0){60}}    \put(75,55){\line(1,0){60}}   \put(75,60){\line(0,-1){5}}
\put(135,60){\line(0,-1){60}}  \put(130,60){\line(0,-1){60}} \put(130,0){\line(1,0){5}}
\put(130,55){$\bullet$}
\put(85,30){$U_{n-4}(q)$}
\put(100,-5){(ii)}

\put(140,70){\line(1,0){65}}  \put(140,65){\line(1,0){65}} \put(140,70){\line(0,-1){5}}
\put(205,70){\line(0,-1){65}}  \put(200,70){\line(0,-1){65}} \put(200,5){\line(1,0){5}}
\put(200,65){$\bullet$}
\put(145,65){\line(1,0){55}}  \put(145,60){\line(1,0){55}} \put(145,65){\line(0,-1){5}}
\put(200,65){\line(0,-1){55}}  \put(195,65){\line(0,-1){55}} \put(195,10){\line(1,0){5}}
\put(195,60){$\bullet$}
\put(150,60){\line(1,0){60}}  \put(150,55){\line(1,0){60}} \put(150,60){\line(0,-1){5}}
\put(210,60){\line(0,-1){60}}  \put(205,60){\line(0,-1){60}} \put(205,0){\line(1,0){5}}
\put(205,55){$\bullet$}
\put(160,30){$U_{n-6}(q)$}
\put(170,-5){(iii)}

\put(220,70){\line(1,0){60}}  \put(220,65){\line(1,0){60}} \put(220,70){\line(0,-1){5}}
\put(280,70){\line(0,-1){60}}  \put(275,70){\line(0,-1){60}} \put(275,10){\line(1,0){5}}
\put(275,65){$\bullet$}
\put(225,65){\line(1,0){65}}  \put(225,60){\line(1,0){65}} \put(225,65){\line(0,-1){5}}
\put(290,65){\line(0,-1){65}}  \put(285,65){\line(0,-1){65}} \put(285,0){\line(1,0){5}}
\put(285,60){$\bullet$}
\put(230,60){\line(1,0){55}}  \put(230,55){\line(1,0){55}} \put(230,60){\line(0,-1){5}}
\put(285,60){\line(0,-1){55}}  \put(280,60){\line(0,-1){55}} \put(280,5){\line(1,0){5}}
\put(280,55){$\bullet$}
\put(230,30){$U_{n-6}(q)$}
\put(245,-5){(iv)}

\put(300,70){\line(1,0){60}}  \put(300,65){\line(1,0){60}} \put(300,70){\line(0,-1){5}}
\put(360,70){\line(0,-1){60}}  \put(355,70){\line(0,-1){60}} \put(355,10){\line(1,0){5}}
\put(355,65){$\bullet$}
\put(305,65){\line(1,0){60}}  \put(305,60){\line(1,0){60}} \put(305,65){\line(0,-1){5}}
\put(365,65){\line(0,-1){60}}  \put(360,65){\line(0,-1){60}} \put(360,5){\line(1,0){5}}
\put(360,60){$\bullet$}
\put(310,60){\line(1,0){60}}  \put(310,55){\line(1,0){60}} \put(310,60){\line(0,-1){5}}
\put(370,60){\line(0,-1){60}}  \put(365,60){\line(0,-1){60}} \put(365,0){\line(1,0){5}}
\put(365,55){$\bullet$}
\put(310,30){$U_{n-6}(q)$}
\put(325,-5){(v)}

\end{picture}

\end{center}
The condition of cases (iii), (iv) and (v) is $n\geq 7$, hence, we suppose $n\geq 7$ in this subsection.

In (i), by Theorem \ref{thm:unique_factor}, each tensor product of an almost faithful elementary character with a
third largest degree irreducible character of $U_{\alpha_{2,n-2}}\simeq
U_{n-2}(q)$ is a third largest degree almost faithful irreducible
character of $U$. In this case, there are
$(q-1)N_{n-2,\mu(n-2)-2}(q)$ third largest degree irreducible characters.

In (ii), 
a tensor product of  two elementary characters  associated to $\alpha_{1,n-2}$
and $\alpha_{2,n-1}$ respectively decomposes into $q$ distinct
irreducible constituents. A tensor product of each constituent
with a second largest degree irreducible of
$U_{\alpha_{3,n-3}}\simeq U_{n-4}(q)$ is a third largest degree
irreducible character of $U$. Hence, in this case, they are not almost faithful and there
are $q(q-1)^2N_{n-4,\mu(n-4)-1}(q)$ third largest degree irreducible characters.

In (iii), we observe a tensor product of three elementary characters
associated to $\alpha_{2,n-3}$, $\alpha_{1,n-2}$ and $\alpha_{3,n-1}$ respectively; and in
(iv), a tensor product of three elementary characters associated to
$\alpha_{1,n-3}$, $\alpha_{3,n-2}$ and $\alpha_{2,n-1}$ respectively. Since cases
(iii) and (iv) are symmetric by the graph automorphism of $\GL_n(q)$, their decompositions are similar.


\begin{slemma} \label{slem:case(iii)}

A tensor product of three elementary characters associated to $\alpha_{2,n-3}$,
$\alpha_{1,n-2}$ and $\alpha_{3,n-1}$ respectively, in case (iii), decomposes
into $q^2$ distinct irreducible constituents of degree $q^{3n-14}$, each has
multiplicity $1$. A tensor product of each irreducible constituent with a
largest degree irreducible character of $U_{\alpha_{4,n-4}}\simeq
U_{n-6}(q)$ is a third largest degree irreducible character of $U$.
So there are $q^2(q-1)^3N_{n-6,\mu(n-6)}(q)$ third largest degree irreducible characters in this case.

\end{slemma}

\begin{proof}
Let
$D:=\{\alpha_{2,n-3},\alpha_{1,n-2},\alpha_{3,n-1}\}$, $V_1:=V_{\alpha_{2,n-3}}$,
$V_2:=V_{\alpha_{1,n-2}}$, $V_3:=V_{\alpha_{3,n-1}}$, and
$\phi:=\{\lam_1,\lam_2,\lam_3\}$ where $\lam_i\in
\Irr(V_i/[V_i,V_i])^*$ for $i=1,2,3$. By Lemma \ref{lem:tobelinear}, we have
$\xi_{D,\phi}=\bigotimes_{\lam_i\in \phi}{\lam_i}^U=\lam^U$ where
$\lam=\bigotimes_{\lam_i\in \phi}\lam_i|_{V_D}$ and $V_D=V_1\cap
V_2\cap V_3$.

Let $X:=X_{\alpha_{1,2}}X_{\alpha_2}$. Since $[X,V_D]=\{1\}$, $\la X,V_D\ra
=XV_D\simeq X\times V_D$. Since $X$ is abelian, 
$\lam^{XV_D}$ splits into $q^2$ distinct linear characters. Let $\mu\neq\mu'$ be two extensions of
$\lam$ to $XV_D$.

Let $M_0:=XV_D$, $M_i:=M_{i-1}X_{\alpha_{2,n-3-i}}$ for $i\in[1,n-5]$, $M_i:=M_{i-1}X_{\alpha_{1,2n-7-i}}$ for $i\in[n-4,2n-8]$, and
$M_i:=M_{i-1}X_{\alpha_{3,3n-9-i}}$ for $i\in[2n-7,3n-12]$. So $M_{n-4}=M_{n-5}$, $M_{2n-6}=M_{2n-7}$, $M_{n-5}=(V_2\cap
V_3)X$, $M_{2n-8}=V_3$, $M_{3n-12}=U$.
Let $K_0:=(leg(\alpha_{2,n-3})\cup leg(\alpha_{1,n-2})\cup
leg(\alpha_{3,n-1}))-\{\alpha_{3,n-3},\alpha_{3,n-2}\}$, and
$K_i:=K_{i-1}-\{\alpha_{n-2-i,n-3}\}$ for $i\in[1,n-5]$,
$K_i:=K_{i-1}-\{\alpha_{2n-6-i,n-2}\}$ for $i\in[n-4,2n-8]$,
$K_i:=K_{i-1}-\{\alpha_{3n-8-i,n-1}\}$ for $i\in[2n-7,3n-12]$. It is clear that $K_{n-4}=K_{n-5}$, $K_{2n-6}=K_{2n-7}$,
$K_{3n-14}$ is empty, $X_\beta\subset \ker(\mu)$, $\ker(\mu')$ for all $\beta\in
K_0$; and for every $M_i=M_{i-1}X_\alpha$, $M_i\neq M_{i-1}$, there is
unique $\beta\in K_i- K_{i-1}$ such that $\alpha+\beta\in
\{\alpha_{2,n-3},\alpha_{1,n-1},\alpha_{3,n-3}\}$.

Now we achieve $\mu^U\neq \mu'^U\in \Irr(U)$
by applying exactly the same technique as in Lemma
\ref{lem:prod2midafis} (ii), i.e. by a sequence of inductions along the arms of
$\alpha_{2,n-3}$, $\alpha_{1,n-2}$, $\alpha_{3,n-1}$ respectively, namely
from $M_0$ to $M_1$, ..., $M_{3n-12}$.
\end{proof}

\smallskip
In both cases (iii) and (iv), there are
$2q^2(q-1)^3N_{n-6,\mu(n-6)}(q)$ characters. Next, in (v), we
observe a tensor product of three elementary characters associated to
$\alpha_{1,n-3}$, $\alpha_{2,n-2}$ and $\alpha_{3,n-1}$ respectively.

\begin{slemma} \label{slem:case(v)}

A tensor product of three elementary characters associated to $\alpha_{1,n-3}$,
$\alpha_{2,n-2}$ and $\alpha_{3,n-1}$ respectively, in case (v), decomposes into $q^2$ distinct irreducibles of degree $q^{3n-15}$ with multiplicity $1$  and $(q-1)$ distinct irreducibles of degree $q^{3n-14}$ with multiplicity
$q$. A tensor product of each constituent of degree $q^{3n-14}$ with a largest
degree irreducible character of $U_{\alpha_{4,n-4}}\simeq U_{n-6}(q)$ is a third largest
degree irreducible character of $U$. Hence, there are
$(q-1)^4N_{n-6,\mu(n-6)}(q)$ third largest degree irreducible characters in this case.

\end{slemma}

\begin{proof}
Let
$D:=\{\alpha_{1,n-3},\alpha_{2,n-2},\alpha_{3,n-1}\}$, $V_1:=V_{\alpha_{1,n-3}}$,
$V_2:=V_{\alpha_{2,n-2}}$, $V_3=V_{\alpha_{3,n-1}}$, and
$\phi:=\{\lam_1,\lam_2,\lam_3\}$ where $\lam_i\in
\Irr(V_i/[V_i,V_i])^*$, $i=1,2,3$. By Lemma \ref{lem:tobelinear},
$\xi_{D,\phi}=\lam^U$ where $\lam=\bigotimes_{\lam_i\in
\phi}\lam_i|_{V_D}$ and $V_D=V_1\cap V_2\cap V_3$.

Let $X:=X_{\alpha_1}X_{\alpha_{1,2}}X_{\alpha_2}$. Since $[X,V_D]=\{1\}$, $\la
X,V_D\ra 
\simeq X\times V_D$ and $\lam$ inflates to $XV_D$. We have
$\lam^{XV_D}=(\lam_{XV_D}|_{V_D}\otimes
1_{V_D})^{XV_D}=\lam_{XV_D}\otimes
{1_{V_D}}^{XV_D}=\lam_{XV_D}\otimes \sum_{\chi\in
\Irr(XV_D/V_D)}\chi(1)\cdot \chi$.
Since ${1_{V_D}}^{XV_D}=({1_{\{1\}}}^X)_{XV_D}$, the
regular character of $X$ inflated to $XV_D$ and $X\simeq U_3(q)$,
${1_{V_D}}^{XV_D}$ decomposes into $(q-1)$ irreducibles of degree
$q$ and $q^2$ linear characters. Since $\Irr(XV_D)\simeq \Irr(X)\times
\Irr(V_D)$, $\lam_{XV_D}\otimes \chi_{XV_D}\in \Irr(XV_D)$, for all $\chi\in
\Irr(X)$. Let $\mu\neq \mu'$ be two irreducible constituents of $\lam^{XV_D}$.

Let $M_0:=XV_D$, and $M_i:=M_{i-1}X_{\alpha_{1,n-3-i}}$ for $i\in[1,n-4]$, $M_i:=M_{i-1}X_{\alpha_{2,2n-6-i}}$ for $i\in[n-3,2n-8]$,
$M_i:=M_{i-1}X_{\alpha_{3,3n-9-i}}$ for $i\in[2n-7,3n-12]$. Hence, $M_{n-2}=M_{n-3}=M_{n-4}$, $M_{2n-9}=M_{2n-8}$,
$M_{n-4}=(V_2\cap V_3)X$, $M_{2n-8}=V_3$, and $M_{3n-12}=U$.
Let $K_0:=(leg(\alpha_{1,n-3})\cup leg(\alpha_{2,n-2})\cup
leg(\alpha_{3,n-1}))-\{\alpha_{2,n-3},\alpha_{3,n-3},\alpha_{3,n-2}\}$,
and $K_i:=K_{i-1}-\{\alpha_{n-2-i,n-3}\}$ for $i\in[1,n-4]$,
$K_i:=K_{i-1}-\{\alpha_{2n-5-i,n-2}\}$ for $i\in[n-3,2n-8]$,
$K_i:=K_{i-1}-\{\alpha_{3n-8-i,n-1}\}$ for $i\in[2n-7,3n-12]$. Hence, $K_{n-2}=K_{n-3}=K_{n-4}$,
$K_{2n-9}=K_{2n-8}$, $K_{3n-12}$ is empty, $X_\beta\subset \ker(\mu)$, $\ker(\mu')$ for all $\beta\in K_0$;  for every
$M_i=M_{i-1}X_\alpha$ such that $M_i\neq M_{i-1}$, there is unique $\beta\in K_i- K_{i-1}$ such that $\alpha+\beta\in
\{\alpha_{1,n-3},\alpha_{2,n-2},\alpha_{3,n-1}\}$.

We get $\mu^U\neq \mu'^U\in \Irr(U)$ by applying the same technique as in Lemma \ref{lem:prod2midafis} (ii), i.e. proceeding a sequence of inductions along the arms of $\alpha_{1,n-3}$, $\alpha_{2,n-2}$, $\alpha_{3,n-1}$ respectively,
namely from $M_0$ to $M_1$, $M_2$, ..., $M_{3n-12}$.

Therefore, $\xi_{D,\phi}$ has $(q-1)$ distinct irreducibles of
degree $q^{(3n-12)-3+1}=q^{3n-14}$ with multiplicity $q$ and $q^2$
distinct irreducibles of degree $q^{3n-15}$ with multiplicity
$1$.
\end{proof}

\begin{stheorem}\label{slem:third_largest}

For $n\geq 7$, each third largest degree irreducible representation of $U$ is of the form (i), (ii), (iii), (iv), or (v). Therefore,

\[\begin{array}{ll}
N_{n,\mu(n)-2}(q)=&(q-1)N_{n-2,\mu(n-2)-2}(q)+q(q-1)^2N_{n-4,\mu(n-4)-1}(q)\\
&+(2q^2+q-1)(q-1)^3N_{n-6,\mu(n-6)}(q).
\end{array}\]
\end{stheorem}

\begin{proof}
It is clear that the recursion
formula we want to prove comes from (i), ..., (v), hence we
mainly show that there are no more supercharacters $\xi_{D,\phi}$ which give
irreducible constituents of degree $q^{\mu(n)-2}=q^{n-4+\mu(n-2)}$.

Case (i) lists all third largest degree almost faithful irreducible
characters; cases (ii), ..., (v) correspond with cases where
$X_{\alpha_{1,n-1}}$, $X_{\alpha_{1,n-1}}X_{\alpha_{1,n-2}}$ or $X_{\alpha_{1,n-1}}X_{\alpha_{2,n-1}}$ is contained in the
kernel of $\xi_{D,\phi}$. By Theorem \ref{thm:Andre} and the
graph automorphism of $U$, it suffices to check two cases: the first type $\xi_{D,\phi}$  with $\alpha_{1,n-1-t}\in D$,
$t\geq 3$ (corresponding to cases (iii) and (iv)), i.e. $X_{\alpha_{1,n-1}}X_{\alpha_{1,n-2}}X_{\alpha_{1,n-3}}\subset \ker(\xi_{D,\phi})$, and the second type $\xi_{D,\phi}$ with
$\{\alpha_{1,n-3},\alpha_{3,n-1}\}\subset D$ and $X_{\alpha_{1,n-2}}
X_{\alpha_{1,n-1}}X_{\alpha_{2,n-2}}X_{\alpha_{2,n-1}}\subset
\ker(\xi_{D,\phi})$ (corresponding to case (v)).

First, we study $\xi_{D,\phi}$  with $\alpha_{1,n-1-t}\in D$,  $t\geq 3$. Let $\chi\in \Irr(\xi_{D,\phi})$. Since
$\prod_{i=1}^{3} X_{\alpha_{1,n-i}}\subset \ker(\xi_{D,\phi})$, we
proceed $U$ by modulo $\prod_{i=1}^3 X_{\alpha_{1,n-i}}$.
If $\alpha_{k+1,n-1}\notin D$ for all $k\in[1,n-2]$, $\chi$ can be considered as a character of $U_{\alpha_{1,n-2}}\simeq U_{n-1}(q)$. By Lemma \ref{slem:largest} and Lemma \ref{slem:second_largest}, $\chi$ has degree at most $q^{\mu(n-1)-2}$.
By Lemma \ref{slem:mu_n-1}, $q^{\mu(n-1)-2}<q^{\mu(n)-2}$ for all $n\geq 7$, it is done. So we suppose that $\alpha_{k+1,n-1}\in D$ for some $k\in[1,n-2]$. Set $\alpha:=\alpha_{k+1,n-1}$.

By Lemma \ref{p2}, $\chi=(\xi_{V_\alpha}\otimes \lambda)^U$ for some $\xi\in \Irr(V_\alpha\cap U_{\alpha_{1,n-2}})$ and $\lambda\in \Irr(V_\alpha/[V_\alpha,V_\alpha])^*$. Hence, $\chi(1)=\xi(1)\lambda^U(1)=q^{n-k-2}\xi(1)$. We shall show that $\xi(1)<q^{\mu(n-2)+k-2}$.

Let $T:=\la X_\beta\subset V_\alpha \cap U_{\alpha_{1,n-2}}:\ \beta\notin leg(\alpha_{1,k})\ra = T_{k,2}\times X_{\alpha_{1,k}}\simeq U_{n-2}(q)\times X_{\alpha_{1,k}}$ by Lemma \ref{p1}. Since we work with $U/\prod_{i=1}^{3} X_{\alpha_{1,n-i}}$, by Lemma \ref{slem:largest} and Lemma \ref{slem:second_largest}, $T_{k,2}$ only has irreducible characters of degree at most $q^{\mu(n-2)-2}$. Since $|V_\alpha:T|=q^{|leg(\alpha_{1,k})|}=q^{k-1}$ and each irreducible character of $V_\alpha$ is a constituent of some induced character of $T$, $\xi(1)\leq q^{\mu(n-2)-2+k-1}<q^{\mu(n-2)+k-2}$.

Now we study the second type  of $\xi_{D,\phi}$ with
$\{\alpha_{1,n-3},\alpha_{3,n-1}\}\subset D$ and $X_\beta\subset
\ker(\xi_{D,\phi})$ for all $\beta\in
S=\{\alpha_{1,n-2},\alpha_{2,n-2},\alpha_{1,n-1},\alpha_{2,n-1}\}$. We proceed $U$ by modulo
$\prod_{\beta\in S}X_\beta$. Let $\chi\in \Irr(\xi_{D,\phi})$ and set $\alpha:=\alpha_{3,n-1}$.

With $k=2$, we repeat the proof with \[T:=\la X_\beta\subset V_\alpha \cap U_{\alpha_{1,n-2}}:\ \beta\notin leg(\alpha_{1,2})\ra = T_{2,2}\times X_{\alpha_{1,2}}\simeq U_{n-2}(q)\times X_{\alpha_{1,2}}.\] Since we work with $U/X_{\alpha_{1,n-2}}X_{\alpha_{2,n-2}}$, by Lemma \ref{slem:largest} and Lemma \ref{slem:second_largest}, $T_{2,2}$ only has irreducible characters of degree at most $q^{\mu(n-2)-2}$. Therefore, the claim holds.
\end{proof}

\smallskip
This method can be generalized to find $N_{n,\mu(n)-4}(q)$ for $n\geq 9$. First, all supercharacters $\xi_{D,\phi}$ with $|D|=4$ must be decomposed. To obtain large degree constituents, $D$ should be chosen from $\{\alpha_{i,j}:1\leq i\leq 4, n-4\leq j\leq n-1\}$. By Lemma \ref{lem:D_decomposable}, one only needs to consider indecomposable basic sets $D$ since the others belong to the largest, second largest and third largest degree representations. It is easy to check that all these characters have degree $q^{4n-20}=q^{\mu(n)-\mu(n-8)}$. One shall obtain $(4q^3+4q^2-3q-1)(q-1)^4$ distinct irreducible constituents of degree $q^{4n-23}$ after decomposing them all. Finally, the technique in the proof of Theorem \ref{slem:third_largest} is used
to prove that there are no more supercharacters giving constituents of degree equal to $q^{\mu(n)-3}$. Hence, we obtain the formula.
\[
\begin{array}{ll}
N_{n,\mu(n)-3}(q)=&(q-1)N_{n-2,\mu(n-2)-3}(q)+q(q-1)^2N_{n-4,\mu(n-4)-2}(q)\\
&+(2q^2+q-1)(q-1)^3N_{n-6,\mu(n-6)-1}(q)+q^2(q-1)^3 N_{n-6,\mu(n-6)}(q)\\
&+(4q^3+4q^2-3q-1)(q-1)^4N_{n-8,\mu(n-8)}(q).
\end{array}\]

\section*{Acknowledgement}
The author would like to thank Professor
Kay Magaard for all his help to approach
complex representations of Sylow $p$-subgroups of finite
groups of Lie type $G(p^f)$ by their root systems. This paper
was presented at the Glauberman Conference 2008 at the University of
Chicago to dedicate to Glauberman's birthday. The author would like to thank the referee for many helpful comments and suggestions.





\bibliographystyle{elsarticle-num}







\end{document}